\documentclass[preprint,3p]{elsarticle}

\usepackage{amssymb}
\usepackage{amsmath}
\usepackage{amsthm}
\usepackage{mathtools}
\usepackage{mathrsfs}
\usepackage{enumitem}
\usepackage[colorlinks=true]{hyperref}
\usepackage[nameinlink,capitalise,noabbrev]{cleveref}

\makeatletter
\renewcommand\subsection{\@startsection{subsection}{2}{\z@}%
  {15\p@ \@plus 6\p@ \@minus 3\p@}%
  {6\p@ \@plus 3\p@ \@minus 2\p@}%
  {\normalfont\normalsize\bfseries\boldmath}}
\renewcommand\subsubsection{\@startsection{subsubsection}{3}{\z@}%
  {12\p@ \@plus 6\p@ \@minus 3\p@}%
  {4\p@ \@plus 2\p@ \@minus\p@}%
  {\normalfont\normalsize\bfseries\boldmath}}
\makeatother

\setlist[enumerate]{font=\upshape,label=(\roman*)}
\numberwithin{equation}{section}

\theoremstyle{plain}
\newtheorem{theorem}{Theorem}[section]
\newtheorem{lemma}[theorem]{Lemma}
\newtheorem{proposition}[theorem]{Proposition}
\newtheorem{corollary}[theorem]{Corollary}
\theoremstyle{definition}

\newtheorem{remark}[theorem]{Remark}

\newcommand{\N}{\mathbb N}
\newcommand{\Nzero}{\mathbb N_0}
\newcommand{\Z}{\mathbb Z}
\newcommand{\Q}{\mathbb Q}
\newcommand{\PP}{\mathbb P}
\newcommand{\De}{\mathcal D_e}
\newcommand{\taue}{\tau_e}
\newcommand{\Homeo}{\operatorname{Homeo}}
\newcommand{\Sym}{\operatorname{Sym}}
\newcommand{\lcm}{\operatorname{lcm}}
\newcommand{\rad}{\operatorname{rad}}
\newcommand{\ord}{\operatorname{ord}}

\begin{document}

\begin{frontmatter}

\title{Rigidity and homeomorphism groups of arithmetic progression topologies}

\author{Yifan Meng}
\ead{mengyifanmoe@qq.com}

\address{School of Mathematics, Sichuan University, Chengdu 610064, China}

\begin{abstract}
The Golomb and Kirch topologies on the positive integers are both rigid, meaning that the identity is their only self-homeomorphism. We interpolate between them using the subbasis of arithmetic progressions $a+p^k\Nzero$, where $p$ is prime, $a$ is a positive integer with $p\nmid a$, and $k\in\Nzero$ is bounded by a height function $e(p)\in\Nzero\cup\{\infty\}$. The constant functions $e\equiv1$ and $e\equiv\infty$ recover the Kirch and Golomb topologies, respectively. We prove that the resulting space is rigid whenever $e(p)\ge1$ for all but finitely many primes. When only finitely many primes have positive height, we determine the full self-homeomorphism group and establish necessary and sufficient conditions for rigidity. We also construct a nonrigid space with $e=1$ on a set of natural density one among the primes and $e=0$ elsewhere.
\end{abstract}


\end{frontmatter}

\section{Introduction}
\label{sec:introduction}

Banakh, Spirito, and Turek proved that the Golomb topology on the positive integers is rigid~\cite{Banakh2021a}, building on earlier work of Banakh, Mioduszewski, and Turek~\cite{Banakh2018}. Shortly thereafter, Banakh, Stelmakh, and Turek established rigidity for the Kirch topology as well~\cite{Banakh2021b}. In both spaces, the identity is the sole self-homeomorphism.

Write $\N=\{1,2,\ldots\}$ and $\Nzero=\N\cup\{0\}$. Both topologies are generated by subbasic arithmetic progressions $a+p^k\Nzero$ with $a\in\N$, $p$ prime, and $p\nmid a$; the Golomb topology allows all $k\in\Nzero$~\cite{Golomb1959}, whereas the Kirch topology restricts to $k\in\{0,1\}$~\cite{Kirch1969}.

Let $\PP$ denote the set of all primes. To interpolate between these two topologies, we specify an exponent bound at each prime via a \emph{height function}
\[
 e\colon\PP\longrightarrow\Nzero\cup\{\infty\}.
\]
The function $e$ defines a topology $\taue$ on $\N$ with subbasis
\begin{equation}\label{eq:prime-subbase}
 a+p^k\Nzero,\qquad
 a\in\N,\quad p\in\PP,\quad k\in\Nzero,\quad k\le e(p),\quad p\nmid a.
\end{equation}
We write $X_e=(\N,\taue)$. The constant functions $e\equiv1$ and $e\equiv\infty$ yield the Kirch topology $\tau_{\mathrm K}$ and the Golomb topology $\tau_{\mathrm G}$, respectively. Whenever $e(p)\ge1$ for every prime $p$, we have
\[
 \tau_{\mathrm K}\subseteq\taue\subseteq\tau_{\mathrm G}.
\]

In the broader framework of coset topologies on $\Z$ introduced by Longhi, Mu, and Saettone~\cite[Definition~2.28, Example~2.29, and Remark~2.30(6)]{Longhi2023}, the cosets $a+p^k\Z$ with $1\le k\le e(p)$ and $p\nmid a$ form a subbasis. In their notation, this corresponds to setting the upper bound $\kappa=e$, taking $\gamma(p)=1$ and $B_p=\{0\}\subseteq\Z/p\Z$ at primes with $e(p)>0$ (which starts at $k=1$ and excludes multiples of $p$), and taking $\gamma(p)=0$ at primes with $e(p)=0$ (which omits all cosets at $p$). Restricting these cosets to $\N$ produces the corresponding unit residue classes. When infinitely many primes have positive height, these classes generate $\taue$ (\cref{prop:inf-terminal-basis}); when only finitely many do, the upward rays $a+\Nzero$ make $\taue$ strictly finer than the restricted coset topology.

Special cases of this family have been studied previously: Szyszkowska investigated constant finite height functions $e\equiv m$~\cite{Szyszkowska2022a,Szyszkowska2022b}, and Kwela and Tryba generalized these constructions using arbitrary sequences of nonnegative integers to compare their ideals of nowhere dense sets~\cite{KwelaTryba2023}.

Define the \emph{support} of $e$ by
\begin{equation}\label{eq:height-partition-global}
 S_e=\{p\in\PP:e(p)>0\}.
\end{equation}
We call the primes in $S_e$ \emph{active} and address three natural questions:
\begin{enumerate}
\item Is $X_e$ rigid whenever every prime is active?
\item If only finitely many primes are active, when is $X_e$ rigid, and what is its self-homeomorphism group?
\item Does infinite support alone guarantee rigidity?
\end{enumerate}

We answer the first question affirmatively in a stronger form: $X_e$ is rigid whenever $e(p)\ge1$ for all but finitely many primes (\cref{thm:main-cofinite-arbitrary}), and in particular whenever $e(p)\ge1$ for every prime $p$ (\cref{cor:full-support-interpolation}). For finite support, we completely determine $\Homeo(X_e)$ and obtain necessary and sufficient conditions for rigidity (\cref{sec:finite-support}). We answer the third question in the negative even for a support of natural density one among the primes. More precisely, we construct a coinfinite set $G\subseteq\PP$ with $|G\cap[1,x]|/\pi(x)\to1$, where $\pi$ is the prime-counting function, and a fixed-point-free involution of $\N$ exchanging $1$ and $2$. For every infinite $S\subseteq G$, this involution is a self-homeomorphism of the space with $e=1$ on $S$ and $e=0$ elsewhere (\cref{thm:sparse-involution}).

Our proof of cofinite rigidity builds on the filter method of Banakh, Stelmakh, and Turek~\cite[\S2 and Lemma~14]{Banakh2021b}, which uses filters generated by neighborhood closures. In the Kirch topology, each prime $p>3$ is identified by counting the $p-1$ filters corresponding to the nonzero residue classes modulo $p$, while $p=3$ is distinguished separately via filter inclusions. When primes are allowed different exponent bounds, however, distinct primes can share the same count of unit residue classes---for instance, there are six unit classes both modulo $9=3^2$ and modulo $7$. Starting from the same filters, we show that every self-homeomorphism permutes the sets $p\N$ for active primes $p$ and transports the unit classes modulo $p^{e(p)}$ at finite height to those at the image prime. When only finitely many primes are inactive, we rule out exchanges between primes with equal totient counts by choosing a third active prime with prescribed residues modulo the two primes (\cref{prop:cof-finite-coordinates-fixed}).

\section{Preliminaries}
\label{sec:framework}

We begin by recording two standard number-theoretic facts that will be used throughout the paper.

\begin{lemma}\label{lem:general-crt}
For $a_i,b_i\in\N$ ($1\le i\le r$, $r\ge1$), the intersection $\bigcap_{i=1}^r(a_i+b_i\Nzero)$ is nonempty if and only if $\gcd(b_i,b_j)\mid a_i-a_j$ for all $i,j$. When nonempty, it equals $a+b\Nzero$, where $b=\lcm(b_1,\ldots,b_r)$ and $a$ is the least element of the intersection.
\end{lemma}
\begin{proof}
If $b=1$, the intersection is $\max_i a_i+\Nzero$. Suppose $b>1$. Any common solution $x$ of $x\equiv a_i\pmod{b_i}$ necessarily satisfies $\gcd(b_i,b_j)\mid a_i-a_j$ for all $i,j$.

Conversely, assume these divisibility conditions hold. For each prime $p\mid b$, choose an index $i$ such that $b_i$ is divisible by the highest power $p^m$ of $p$ dividing $b$. If $p^s\mid b_j$, then $p^s\mid\gcd(b_i,b_j)$, so $a_i\equiv a_j\pmod{p^s}$; hence $x\equiv a_i\pmod{p^m}$ implies $x\equiv a_j\pmod{p^s}$. Because these prime-power moduli are pairwise coprime and multiply to $b$, the Chinese remainder theorem~\cite[Chapter~II, Theorem~2.1]{Lang2002} yields a unique residue class modulo $b$. The intersection consists of all solutions $x\ge\max_i a_i$, and letting $a$ be the least such solution gives $a+b\Nzero$.
\end{proof}

\begin{lemma}\label{lem:dirichlet}
For $a,d\in\N$, the progression $a+d\Nzero$ contains infinitely many primes if and only if $\gcd(a,d)=1$.
\end{lemma}
\begin{proof}
If $\gcd(a,d)=1$, Dirichlet's theorem~\cite[Chaps.~1 and~4]{Davenport2000} guarantees infinitely many primes congruent to $a$ modulo $d$, and only finitely many of them lie below $a$. If $\gcd(a,d)>1$, every term of $a+d\Nzero$ is divisible by $\gcd(a,d)$, so at most one term can be prime.
\end{proof}

Let $v_p(d)$ denote the $p$-adic valuation of $d\in\N$. Given a height function $e$, define the set of \emph{admissible common differences} by
\begin{equation}\label{eq:legal-moduli-global}
 \De=\{d\in\N:v_p(d)\le e(p)\text{ for every }p\in\PP\}.
\end{equation}

\begin{lemma}\label{lem:tail-basis}
The arithmetic progressions
\begin{equation}\label{eq:basic-tail-global}
 a+d\Nzero,\qquad a\in\N,\quad d\in\De,\quad\gcd(a,d)=1,
\end{equation}
form a basis for $\taue$. In particular, those starting at $x$ form a neighborhood basis at $x$.
\end{lemma}
\begin{proof}
The family covers $\N$, and by \cref{lem:general-crt} any nonempty finite intersection of its members is again of the form $a+d\Nzero$, where $d$ is the least common multiple of the original common differences. Thus $d\in\De$ and $\gcd(a,d)=1$. Every subbasic set in~\eqref{eq:prime-subbase} clearly belongs to this family. Conversely, for $d>1$ we have $a+d\Nzero=\bigcap_{p\mid d}(a+p^{v_p(d)}\Nzero)$, while for $d=1$ the upward ray $a+\Nzero$ belongs to the subbasis~\eqref{eq:prime-subbase} upon taking $k=0$ and any prime $p>a$. Finally, if $x\in a+d\Nzero$, then $x+d\Nzero\subseteq a+d\Nzero$ with $\gcd(x,d)=1$, giving the stated neighborhood basis.
\end{proof}

Unless stated otherwise, $\overline A$ denotes the closure of $A$ in $X_e$. We define the specialization preorder on $X_e$ by
\begin{equation}\label{eq:specialization-convention}
 x\preceq_e y\quad\Longleftrightarrow\quad
 x\in\overline{\{y\}}.
\end{equation}
Equivalently, $x\preceq_e y$ if and only if every neighborhood of $x$ contains $y$. Write $\mathord{\uparrow}_e x=\{y\in\N:x\preceq_e y\}$ for the principal upper set of $x$. Recall that a topological space is \emph{Alexandrov} if arbitrary intersections of open sets are open, or equivalently, if every principal upper set $\mathord{\uparrow}_e x$ is open; in an Alexandrov space, the open sets are precisely the upper sets of the specialization preorder~\cite[\S1]{Alexandroff1937}.

\begin{proposition}\label{prop:general-specialization}
Every space $X_e$ is $T_0$. Moreover, $X_e$ is Hausdorff if and only if $S_e$ is infinite, and Alexandrov if and only if $S_e$ is finite and $e(p)<\infty$ for every $p\in S_e$.
\end{proposition}
\begin{proof}
Since every upward ray $x+\Nzero$ is open, $x\preceq_e y$ implies $x\le y$, so $\preceq_e$ is antisymmetric and $X_e$ is $T_0$.

If $S_e$ is finite, set $P=\prod_{p\in S_e}p$. The only admissible common difference coprime to $P$ is $1$, so every neighborhood of $P$ contains $P+1$, and $X_e$ is not Hausdorff. If in addition $e(p)<\infty$ for all $p\in S_e$, then every admissible common difference coprime to $x\in\N$ divides $\prod_{p\in S_e,\,p\nmid x}p^{e(p)}$, so
\[
 x+\left(\prod_{\substack{p\in S_e\\p\nmid x}}p^{e(p)}\right)\Nzero
\]
is the smallest neighborhood of $x$, making $X_e$ Alexandrov. Conversely, if $e(p)=\infty$ for some $p$, the neighborhoods $1+p^k\Nzero$ ($k\ge1$) intersect in $\{1\}$, which is not open because every basic neighborhood is infinite; hence $X_e$ is not Alexandrov.

If $S_e$ is infinite, any distinct $x,y\in\N$ are separated by the disjoint neighborhoods $x+p\Nzero$ and $y+p\Nzero$ for any active prime $p>\max\{x,y\}$, as in Golomb's separation argument~\cite[Theorem~2]{Golomb1962}. In this case every principal upper set is a singleton, which is not open, so $X_e$ is Hausdorff and not Alexandrov.
\end{proof}

For the closure formulas in the Golomb and Kirch topologies, see~\cite[Lemma~2.2]{Banakh2018} and~\cite[Lemma~1]{Banakh2021b}; for the constant finite-height case, see~\cite[Theorem~3.4]{Szyszkowska2022b}; and for a quotient-space formula in coset topologies, see~\cite[Lemma~2.23]{Longhi2023}.

\begin{proposition}\label{prop:closure}
For $a\in\N$ and $d\in\De$ with $\gcd(a,d)=1$,
\begin{equation}\label{eq:closure}
 \overline{a+d\Nzero}
 =\N\cap\bigcap_{p\mid d}
 \left(p\N\cup(a+p^{v_p(d)}\Z)\right),
\end{equation}
with the empty intersection interpreted as $\N$ when $d=1$.
\end{proposition}
\begin{proof}
For the forward inclusion, let $x\in\overline{a+d\Nzero}$. If $p\mid d$ and $p\nmid x$, then $x+p^{v_p(d)}\Nzero$ is a basic neighborhood of $x$ and therefore meets $a+d\Nzero$, forcing $x\equiv a\pmod{p^{v_p(d)}}$.

Conversely, let $x$ belong to the right-hand side of~\eqref{eq:closure}, and let $x+c\Nzero$ be any basic neighborhood of $x$. For every prime $p\mid\gcd(c,d)$ we have $p\nmid x$, so $p^{v_p(d)}\mid x-a$. Since $p^{\min\{v_p(c),v_p(d)\}}\mid p^{v_p(d)}$, it follows that $\gcd(c,d)\mid x-a$, and \cref{lem:general-crt} ensures that $(x+c\Nzero)\cap(a+d\Nzero)\ne\varnothing$.
\end{proof}

A space is \emph{superconnected} if the closures of every finite family of nonempty open sets have nonempty intersection~\cite[\S2 and Corollary~2.4]{Banakh2018}. In the Golomb topology, the common-multiple argument for two closures dates back to Brown~\cite{Brown1953}.

\begin{corollary}\label{cor:superconnected}
Every space $X_e$ is superconnected.
\end{corollary}
\begin{proof}
Let $U_1,\ldots,U_r$ be nonempty open sets, and choose a basic progression $a_i+d_i\Nzero\subseteq U_i$ for each $i$. By~\eqref{eq:closure}, $d_i\N\subseteq\overline{U_i}$ for each $i$, so $\prod_{i=1}^r d_i\in\bigcap_{i=1}^r\overline{U_i}$.
\end{proof}

For an active prime $p\in S_e$, an integer $1\le k\le e(p)$, and $a\in\N$ with $p\nmid a$, define the \emph{cylinder}
\begin{equation}\label{eq:inf-cylinder}
 C_e(a,p,k)=\N\cap(a+p^k\Z).
\end{equation}
For $u\in(\Z/p^k\Z)^\times$, we also write $C_e(u,p,k)$ for $C_e(a,p,k)$ with any positive representative $a$ of $u$. The cosets $a+p^k\Z$ are subbasic in~\cite[Definition~2.28]{Longhi2023}. When $S_e$ is finite, these cylinders generate a strictly coarser topology than $\taue$: each cylinder is a union of basic progressions, whereas any common multiple $b>1$ of the active primes lies in no cylinder (so its only neighborhood in the cylinder topology is $\N$), even though $b+\Nzero$ is a proper neighborhood of $b$ in $\taue$.

\section{Finite prime support}\label{sec:finite-support}

For any height function $e$, let $F_e=\{p\in S_e:e(p)<\infty\}$ and $I_e=\{p\in S_e:e(p)=\infty\}$ be the sets of active primes of finite and infinite height, respectively. Throughout this section, assume that $S_e$ is finite. Empty products are understood to equal $1$.

To describe the principal upper sets in $X_e$, define $Q_e\colon\N\to\N$ by
\begin{equation}\label{eq:Qe-global}
 Q_e(n)=\prod_{\substack{q\in F_e\\q\nmid n}}q^{e(q)}
 \qquad(n\in\N).
\end{equation}

\begin{proposition}
\label{prop:finite-support-specialization}
Suppose that $S_e$ is finite, and put $J_e=\prod_{p\in I_e}p$. For every $n\in\N$,
\begin{equation}\label{eq:finite-support-upper}
 \mathord{\uparrow}_e n=
 \begin{cases}
  \{n\},&J_e\nmid n,\\
  n+Q_e(n)\Nzero,&J_e\mid n.
 \end{cases}
\end{equation}
If $J_e\mid n$, this principal upper set is the smallest neighborhood of $n$; if $J_e\nmid n$, then $n$ has no smallest neighborhood.
\end{proposition}
\begin{proof}
If $J_e\nmid n$, there exists $p\in I_e$ with $p\nmid n$. The neighborhoods $n+p^k\Nzero$ ($k\ge1$) intersect in $\{n\}$, so $\mathord{\uparrow}_e n=\{n\}$, which is not open; hence $n$ has no smallest neighborhood. If $J_e\mid n$, every admissible common difference coprime to $n$ divides $Q_e(n)$, and $n+Q_e(n)\Nzero$ is itself an admissible neighborhood of $n$, so it is the smallest neighborhood by \cref{lem:tail-basis}.
\end{proof}

For a block length $P\in\N$ and an index $j\in\Nzero$, define the $P$-block
\begin{equation}\label{eq:common-P-block}
 B_j(P)=\{jP+1,\ldots,(j+1)P\}.
\end{equation}
Write $P_e=\prod_{p\in S_e}p$.

\begin{proposition}
\label{prop:common-finite-backbone}
Suppose that $S_e$ is finite. Every self-homeomorphism $h$ of $X_e$ satisfies
\begin{equation}\label{eq:common-backbone-fixed}
 h(kP_e)=kP_e\qquad(k\in\N)
\end{equation}
and
\[
 h(B_j(P_e))=B_j(P_e)\qquad(j\in\Nzero).
\]
In particular, $|h(n)-n|<P_e$ for every $n\in\N$.
\end{proposition}
\begin{proof}
By \cref{prop:finite-support-specialization}, $\mathord{\uparrow}_e n$ is cofinite if and only if $P_e\mid n$, in which case $\mathord{\uparrow}_e n=n+\Nzero$ has complement of size $n-1$. Since $h[\mathord{\uparrow}_e n]=\mathord{\uparrow}_e h(n)$, comparing complement sizes forces $h(n)=n$ for every multiple $n$ of $P_e$.

Because both the ray $kP_e+\Nzero$ and its starting point $kP_e$ are invariant under $h$, each initial segment $\{1,\ldots,kP_e\}$ is invariant for $k\in\N$. Taking set differences of successive initial segments shows that every $P_e$-block $B_j(P_e)$ is invariant, so $|h(n)-n|<P_e$ for all $n\in\N$.
\end{proof}

\begin{lemma}\label{lem:fm-blocks}
Put $J_e=\prod_{p\in I_e}p$ and $\mathcal M_e=J_e\N$. Then $\mathcal M_e$ consists precisely of the non-maximal points of the specialization order $\preceq_e$. Every $h\in\Homeo(X_e)$ satisfies $h(\mathcal M_e)=\mathcal M_e$ and
\begin{equation}\label{eq:fm-bounds}
 Q_e(h(n))=Q_e(n)\qquad(n\in\mathcal M_e).
\end{equation}
\end{lemma}
\begin{proof}
By \cref{prop:finite-support-specialization}, $\mathcal M_e$ is the set of points with a nontrivial principal upper set, i.e., the non-maximal points of $\preceq_e$, so $h(\mathcal M_e)=\mathcal M_e$. Fix $n\in\mathcal M_e$ and set
\[
 D=J_e\prod_{q\in F_e}q^{e(q)}.
\]
Choose $j\in\Nzero$ with $jD\ge\max\{n,h(n)\}$ and set $E_j=\{jD+1,\ldots,(j+1)D\}$. Since $P_e\mid D$, $E_j$ is a union of invariant $P_e$-blocks, so $h(E_j)=E_j$. Moreover, both $Q_e(n)$ and $Q_e(h(n))$ divide $D$, and $\min E_j>\max\{n,h(n)\}$. Hence each of the upper sets $n+Q_e(n)\Nzero$ and $h(n)+Q_e(h(n))\Nzero$ meets $E_j$ in exactly one point per period of its common difference. Since $h$ maps one upper set onto the other while preserving $E_j$, we obtain
\[
 \frac{D}{Q_e(n)}
 =|(n+Q_e(n)\Nzero)\cap E_j|
 =|(h(n)+Q_e(h(n))\Nzero)\cap E_j|
 =\frac{D}{Q_e(h(n))},
\]
proving~\eqref{eq:fm-bounds}.
\end{proof}

We will frequently pass to subspaces of multiples $A\N$ by deleting prime divisors of $A$ from the support. For $A\in\N$, define the reduced height function
\begin{equation}\label{eq:deleted-profile-common}
 e^{\langle A\rangle}(p)=
 \begin{cases}
  0,&p\mid A,\\
  e(p),&p\nmid A.
 \end{cases}
\end{equation}

\begin{lemma}
\label{lem:deletion-scaling}
For every $A\in\N$, the scaling map
\[
 \lambda_A\colon X_{e^{\langle A\rangle}}\longrightarrow A\N,
 \qquad \lambda_A(n)=An,
\]
is a homeomorphism onto the subspace $A\N\subseteq X_e$.
\end{lemma}
\begin{proof}
The admissible common differences for basic neighborhoods of $n$ in $X_{e^{\langle A\rangle}}$ are the $d\in\De$ with $\gcd(d,An)=1$. For each such $d$, since $\gcd(A,d)=1$, we have
\begin{equation}\label{eq:deletion-scaling-trace}
 (An+d\Nzero)\cap A\N=A(n+d\Nzero).
\end{equation}
Thus $\lambda_A$ identifies the neighborhood basis at $n$ in $X_{e^{\langle A\rangle}}$ with the subspace neighborhood basis at $An$ in $A\N\subseteq X_e$.
\end{proof}

\subsection{Finite support with finite heights}\label{sec:finite-heights}

In this subsection, suppose that $I_e=\varnothing$. Write $S=S_e$ and $Q=Q_e$, and set
\begin{equation}\label{eq:finite-PM}
 P=\prod_{q\in S}q,\qquad M=\prod_{q\in S}q^{e(q)}.
\end{equation}
By \cref{prop:general-specialization}, $X_e$ is Alexandrov, with specialization order
\begin{equation}\label{eq:finite-order-simple}
 n\preceq_e m\quad\Longleftrightarrow\quad
 m\ge n\ \text{and}\ Q(n)\mid m-n.
\end{equation}
Consequently, a bijection $h\colon\N\to\N$ is a self-homeomorphism of $X_e$ if and only if it is an order automorphism of $(\N,\preceq_e)$.

\begin{proposition}\label{prop:finite-invariants}
Let $h\in\Homeo(X_e)$. Then, for every $n\in\N$,
\[
 Q(h(n))=Q(n),\qquad
 q\mid h(n)\ \Longleftrightarrow\ q\mid n\quad(q\in S).
\]
\end{proposition}
\begin{proof}
Since $\mathcal M_e=\N$, \cref{lem:fm-blocks} gives $Q(h(n))=Q(n)$. For each $q\in S$, we have $q\nmid n$ if and only if $q\mid Q(n)$, which yields the divisibility equivalence.
\end{proof}

\begin{corollary}\label{cor:finite-deletion}
Let $h\in\Homeo(X_e)$. For every positive divisor $A$ of $P$, we have $h(A\N)=A\N$, and the map
\[
 X_{e^{\langle A\rangle}}\longrightarrow X_{e^{\langle A\rangle}},
 \qquad n\longmapsto\frac{h(An)}{A},
\]
is a self-homeomorphism.
\end{corollary}
\begin{proof}
By \cref{prop:finite-invariants}, $h(q\N)=q\N$ for each $q\mid A$. Since $A$ is squarefree, $A\N=\bigcap_{q\mid A}q\N$ is $h$-invariant. The stated map is $\lambda_A^{-1}\circ(h|_{A\N})\circ\lambda_A$, which is a self-homeomorphism of $X_{e^{\langle A\rangle}}$ by \cref{lem:deletion-scaling}.
\end{proof}

\begin{lemma}\label{lem:ray-block}
Let $P,d\in\N$, put $g=\gcd(P,d)$, and suppose that $u<v$ lie in the same $P$-block. If $u+d\Nzero$ and $v+d\Nzero$ contain the same number of points in every $P$-block, then $u$ and $v$ lie in the same $g$-block.
\end{lemma}
\begin{proof}
Suppose for contradiction that $u$ and $v$ lie in different $g$-blocks, and let $w$ be the right endpoint of the $g$-block containing $u$. Then $u\le w<v$. Since $g=\gcd(P,d)$ divides $w$, there exists $m\in\Nzero$ such that $md\equiv-w\pmod P$, so $w+md$ is a multiple of $P$. Now $u+md\le w+md<v+md$, so up to the $P$-block endpoint $w+md$, the progression $u+d\Nzero$ has at least $m+1$ terms whereas $v+d\Nzero$ has at most $m$ terms, contradicting the hypothesis of equal counts in every $P$-block.
\end{proof}

For a prime $\ell$, $Y\in\N$, and $c\in\Z/Y\Z$, define
\begin{equation}\label{eq:common-fiber-group}
 \begin{aligned}
 A_c^{\ell,Y}
 &=\{a\in\{1,\ldots,\ell-1\}:a\equiv c\pmod Y\},\\
 G_{\ell,Y}
 &=\prod_{c\in\Z/Y\Z}\Sym(A_c^{\ell,Y}).
 \end{aligned}
\end{equation}
Equivalently,
\[
 G_{\ell,Y}
 =\{\sigma\in\Sym(\{1,\ldots,\ell-1\}):
       \sigma(a)\equiv a\pmod Y\text{ for every }1\le a<\ell\}.
\]
When $S\ne\varnothing$, we single out its largest prime and write
\begin{equation}\label{eq:pYRL}
 p=\max S,\quad
 Y=\prod_{q\in S\setminus\{p\}}q^{e(q)},\quad
 R=\prod_{q\in S\setminus\{p\}}q,\quad L=p^{e(p)-1}.
\end{equation}

We prove the normal form by strong induction on $|S|$, taking the largest active prime $p=\max S$ as the distinguished coordinate. The key step is to show that a self-homeomorphism is uniquely determined by its restriction to $R\N$; by \cref{cor:finite-deletion}, scaling this subspace by $1/R$ reduces the problem to the one-prime space with support $\{p\}$.

\begin{theorem}\label{thm:finite-normal-form}
Suppose $S\ne\varnothing$. Every self-homeomorphism $h$ of $X_e$ has a unique representation
\begin{equation}\label{eq:normal-map}
 \begin{aligned}
 h(pK)&=pK &&(K\ge1),\\
 h(pK+r)&=pK+\sigma_{[K]}(r)
       &&(K\ge0,\ 1\le r<p),
 \end{aligned}
\end{equation}
where $[K]$ denotes the residue class of $K$ modulo $L=p^{e(p)-1}$, and $\sigma_c\in G_{p,Y}$ for each $c\in\Z/L\Z$. Conversely, every such family $(\sigma_c)_{c\in\Z/L\Z}$ defines a self-homeomorphism of $X_e$. In particular,
\begin{equation}\label{eq:finite-homeo-group}
 \Homeo(X_e)\cong\prod_{c\in\Z/L\Z}G_{p,Y}.
\end{equation}
\end{theorem}
\begin{proof}
\noindent\textbf{Step 1.} We first verify that any map $H\colon X_e\to X_e$ of the form~\eqref{eq:normal-map} is a self-homeomorphism. It is a bijection because it fixes $p\N$ pointwise and permutes each set $\{pK+1,\ldots,pK+p-1\}$. Since every $\sigma_c$ belongs to $G_{p,Y}$, we have $H(n)\equiv n\pmod Y$ for all $n\in\N$. Divisibility by $p$ is also preserved, so $Q(H(n))=Q(n)$.

Suppose $x\preceq_e y$. If $p\mid x$, then $H(x)=x$ and $Q(x)\mid Y$. Since $H$ fixes the multiples of $p$ and preserves the intervals between consecutive multiples, $H(y)\ge x$. Moreover, $H(y)\equiv y\pmod Y$ and $Q(x)\mid y-x$, so
\[
 H(y)-H(x)=(H(y)-y)+(y-x)\equiv 0\pmod{Q(x)},
\]
giving $H(x)\preceq_e H(y)$.

If $p\nmid x$, write $x=pK+r$ with $1\le r<p$. Since $p^{e(p)}\mid Q(x)$, there exists $t\ge0$ such that $y=x+tp^{e(p)}=p(K+tL)+r$. Because $K+tL\equiv K\pmod L$, the same permutation $\sigma_{[K]}$ acts on both $x$ and $y$, so $H(y)-H(x)=y-x$. Together with $Q(H(x))=Q(x)$, this gives $H(x)\preceq_e H(y)$. Replacing each $\sigma_c$ by $\sigma_c^{-1}$ yields the inverse map $H^{-1}$, which has the same form and therefore also preserves $\preceq_e$. Thus $H\in\Homeo(X_e)$.

\smallskip\noindent\textbf{Step 2.} To show that every $h\in\Homeo(X_e)$ has the form~\eqref{eq:normal-map}, we proceed by strong induction on $|S|$. For each nonempty finite support $S$ and every choice of positive finite heights on $S$, we prove two statements simultaneously:
\[
 \begin{aligned}
 \mathsf A(S)&:\quad\text{every self-homeomorphism of }X_e\text{ has the form~\eqref{eq:normal-map}},\\
 \mathsf B(S)&:\quad
 \text{the restriction map }\Homeo(X_e)\longrightarrow\Homeo(R\N),\quad
 h\longmapsto h|_{R\N},\quad\text{is injective}.
 \end{aligned}
\]
By \cref{prop:finite-invariants}, $h$ preserves $q\N$ for each $q\in S\setminus\{p\}$ and hence preserves their intersection $R\N$. By \cref{cor:finite-deletion}, division by $R$ identifies $R\N$ with $X_{e^{\langle R\rangle}}$, whose support is $\{p\}$ with height $e(p)$. Thus, once $\mathsf B(S)$ is established, proving $\mathsf A(S)$ reduces to constructing a map of the form~\eqref{eq:normal-map} that agrees with $h$ on $R\N$.

For the base case $S=\{p\}$, \cref{prop:common-finite-backbone} shows that $h$ fixes $p\N$ and preserves every $p$-block. Hence we may write $h(pK+r)=pK+\sigma_K(r)$ for $1\le r<p$, where $\sigma_K\in\Sym(\{1,\ldots,p-1\})$. Take $K\le K'$ with $K\equiv K'\pmod L$. For each $1\le r<p$, the difference $p(K'-K)$ is divisible by $pL=p^{e(p)}$, so $pK+r\preceq_e pK'+r$. Their images satisfy the same order relation, and $p\nmid h(pK+r)$, so $p^{e(p)}$ divides $h(pK'+r)-h(pK+r)=p(K'-K)+\sigma_{K'}(r)-\sigma_K(r)$. Reducing modulo $p$ gives
\[
 0\equiv\sigma_{K'}(r)-\sigma_K(r)\pmod p.
\]
Since $\sigma_K(r),\sigma_{K'}(r)\in\{1,\ldots,p-1\}$, they must be equal. This proves $\mathsf A(\{p\})$ (here $Y=1$, so $G_{p,Y}=\Sym(\{1,\ldots,p-1\})$), and $\mathsf B(\{p\})$ holds trivially since $R=1$.

\smallskip\noindent\textbf{Step 3.} Assume that $|S|\ge2$ and that $\mathsf A(T)$ and $\mathsf B(T)$ hold for all nonempty proper subsets $T\subsetneq S$. We first show that $h$ fixes $p\N$ pointwise.

Suppose first that $S=\{q,p\}$ with $q<p$, so $P=pq$. By \cref{prop:common-finite-backbone}, $h$ fixes $P\N$ pointwise. Take $x\in p\N\setminus P\N$. By \cref{prop:finite-invariants}, $h(x)$ is also divisible by $p$ and coprime to $q$, so both $\mathord{\uparrow}_e x$ and $\mathord{\uparrow}_e h(x)$ have common difference $q^{e(q)}$. For each $P$-block $B_j(P)$,
\[
 h\bigl((\mathord{\uparrow}_e x)\cap B_j(P)\bigr)
 =(\mathord{\uparrow}_e h(x))\cap B_j(P),
\]
so the two upper sets contain the same number of points in every $P$-block. Moreover, $x$ and $h(x)$ lie in the same $P$-block. Since $\gcd(P,q^{e(q)})=q$, \cref{lem:ray-block} implies that if $h(x)\ne x$, then $x$ and $h(x)$ lie in the same $q$-block. However, distinct multiples of $p$ differ by at least $p>q$, forcing $h(x)=x$.

Now suppose $|S|\ge3$, and apply the induction hypothesis to $T=S\setminus\{p\}$ with $t=\max T$. By \cref{cor:finite-deletion}, the map $g\colon X_{e^{\langle p\rangle}}\to X_{e^{\langle p\rangle}}$ given by $g(n)=h(pn)/p$ is a self-homeomorphism. By $\mathsf A(T)$, $g$ fixes $t\N$ pointwise and satisfies
\[
 g(tk+r)=tk+\alpha_{[k]}(r),\qquad 1\le r<t,
\]
where $[k]$ denotes the residue class of $k$ modulo $t^{e(t)-1}$. Take $x\in p\N$ and write $x=p(tk+r)$ with $0\le r<t$. If $r=0$, then $g(x/p)=x/p$, so $h(x)=x$. If $r>0$, choose $q\in T\setminus\{t\}$. Since $t$ is invertible modulo $q$ and $\gcd(q,t^{e(t)-1})=1$, the Chinese remainder theorem supplies $k'\ge0$ such that
\[
 k'\equiv k\pmod{t^{e(t)-1}},\qquad tk'+r\equiv0\pmod q.
\]
Set $y=p(tk'+r)$. Then $q\mid y$ and $[k']=[k]$, so the formula for $g$ gives
\[
 h(x)-x=p\bigl(\alpha_{[k]}(r)-r\bigr)=h(y)-y.
\]
By \cref{cor:finite-deletion}, the map $n\mapsto h(qn)/q$ is a self-homeomorphism of $X_{e^{\langle q\rangle}}$, whose support $S\setminus\{q\}$ still has largest prime $p$. Applying $\mathsf A(S\setminus\{q\})$ shows that this map fixes the multiple $y/q$ of $p$. Hence $h(y)=y$, and therefore $h(x)=x$.

\smallskip\noindent\textbf{Step 4.} Next, we prove that for points coprime to $P$,
\begin{equation}\label{eq:coprime-Y}
 \gcd(z,P)=1\quad\Longrightarrow\quad h(z)\equiv z\pmod Y.
\end{equation}
First assume $z>pY$. Since $\gcd(p,Y)=1$, the Chinese remainder theorem gives $1\le x\le pY$ such that $p\mid x$ and $x\equiv z\pmod Y$. For every $q\in S\setminus\{p\}$, we have $q\nmid z$ and $x\equiv z\pmod q$, so $q\nmid x$. Thus $Q(x)=Y$, and since $x<z$, we have $x\preceq_e z$. As $h(x)=x$ by Step~3, order preservation gives $x\preceq_e h(z)$, so $h(z)\equiv x\equiv z\pmod Y$.

If $z\le pY$ and $\gcd(z,P)=1$, choose $m\in\N$ such that $w=z+mM>pY$. Since $P\mid M$, we have $\gcd(w,P)=1$, and $z\preceq_e w$ because $Q(z)=M$. Order preservation and $Q(h(z))=Q(z)=M$ give $h(w)\equiv h(z)\pmod M$. Since $Y\mid M$ and $h(w)\equiv w\pmod Y$ by the previous paragraph, we conclude that $h(z)\equiv z\pmod Y$.

\smallskip\noindent\textbf{Step 5.} We show that every $y\in R\N$ lies in the same $p$-block as $h(y)$ and satisfies $h(y)\equiv y\pmod Y$. If $y\in P\N$, then $h(y)=y$ by \cref{prop:common-finite-backbone}, so both assertions hold trivially. Assume now that $y\in R\N\setminus P\N$, and write $y=R(pk+r)$ with $k\ge0$ and $1\le r<p$. Both $y$ and $h(y)$ lie in $B_k(P)$, and by \cref{prop:finite-invariants} both are divisible by $R$ and coprime to $p$, so their principal upper sets have common difference $p^{e(p)}$. Since $h$ maps $\mathord{\uparrow}_e y$ onto $\mathord{\uparrow}_e h(y)$ while preserving every $P$-block, the two upper sets have equal counts in every $P$-block. As $\gcd(P,p^{e(p)})=p$, \cref{lem:ray-block} places $y$ and $h(y)$ in the same $p$-block.

To prove $h(y)\equiv y\pmod Y$, set $z=y+p^{e(p)}$. This point is coprime to $P$: since $p\nmid y$, we have $p\nmid z$, and for every $q\mid R$, $z\equiv p^{e(p)}\not\equiv0\pmod q$. Because $Q(y)=p^{e(p)}$ and $y\preceq_e z$, injectivity and order preservation give
\[
 h(z)=h(y)+\ell p^{e(p)}\qquad(\ell\ge1).
\]
By~\eqref{eq:coprime-Y}, $Y\mid h(z)-z$. Since $R\mid Y$ and $R\mid h(y)-y$, it follows that $R\mid(\ell-1)p^{e(p)}$. As $\gcd(R,p)=1$, we may write $\ell=1+nR$ for some integer $n\ge0$.

Set $z_0=h(y)+p^{e(p)}$. Translating the common $p$-block of $y$ and $h(y)$ by $p^{e(p)}$ places $z$ and $z_0$ in the same $p$-block, and hence in the same $P$-block. On the other hand,
\[
 h(z)=z_0+nRp^{e(p)}=z_0+nPL,
\]
so shifting $z_0$ by $nPL$ changes its $P$-block index by $nL$. Since $h(z)$ and $z$ lie in the same $P$-block, we must have $n=0$. Thus $h(y)-y=h(z)-z\equiv0\pmod Y$.

\smallskip\noindent\textbf{Step 6.} We now prove $\mathsf B(S)$, showing that the restriction $h|_{R\N}$ determines $h$ uniquely. If $h_1,h_2\in\Homeo(X_e)$ agree on $R\N$, then $h_2^{-1}\circ h_1$ fixes $R\N$ pointwise; hence it suffices to show that any $h\in\Homeo(X_e)$ fixing $R\N$ pointwise is the identity.

For each $q\in S\setminus\{p\}$, the induced self-homeomorphism $n\mapsto h(qn)/q$ of $X_{e^{\langle q\rangle}}$ fixes $(R/q)\N$ pointwise. Since $R/q$ is the product of the primes other than $p$ in $S\setminus\{q\}$, the inductive hypothesis $\mathsf B(S\setminus\{q\})$ implies that $h$ fixes $q\N$ pointwise. We already know from Step~3 that $h$ fixes $p\N$ pointwise.

It remains to show that $h$ fixes every integer $z$ coprime to $P$. Suppose first that $z>Rp^{e(p)}$. Since $\gcd(R,p)=1$, the Chinese remainder theorem gives $1\le y\le Rp^{e(p)}$ with
\[
 y\equiv0\pmod R,\qquad y\equiv z\pmod{p^{e(p)}}.
\]
Then $y\in R\N$, so $h(y)=y$, and $p\nmid y$, so $Q(y)=p^{e(p)}$. Since $y<z$, we have $y\preceq_e z$ and hence $y\preceq_e h(z)$, which gives $h(z)\equiv y\equiv z\pmod{p^{e(p)}}$. Combined with $h(z)\equiv z\pmod Y$ from Step~4 and $\gcd(Y,p)=1$, this yields $M\mid h(z)-z$. By \cref{prop:common-finite-backbone}, $|h(z)-z|<P\le M$, forcing $h(z)=z$.

If $z\le Rp^{e(p)}$ with $\gcd(z,P)=1$, choose $m\ge1$ such that $w=z+mM>Rp^{e(p)}$. Since $P\mid M$, $w$ is coprime to $P$ and therefore fixed by $h$. From $z\preceq_e w$ we obtain $h(z)\preceq_e w$, and $Q(h(z))=M$ gives $h(z)\equiv w\equiv z\pmod M$. Again $|h(z)-z|<P\le M$ forces $h(z)=z$, establishing $\mathsf B(S)$.

\smallskip\noindent\textbf{Step 7.} Finally, we construct a map of the form~\eqref{eq:normal-map} agreeing with $h$ on $R\N$. By \cref{cor:finite-deletion}, the map
\[
 g_R\colon X_{e^{\langle R\rangle}}\longrightarrow X_{e^{\langle R\rangle}},
 \qquad g_R(n)=\frac{h(Rn)}{R},
\]
is a self-homeomorphism of the one-prime space with support $\{p\}$. By the base case, there exist permutations $\beta_c\in\Sym(\{1,\ldots,p-1\})$ ($c\in\Z/L\Z$) such that
\begin{equation}\label{eq:deep-layer}
 g_R(pk+r)=pk+\beta_{[k]}(r)
 \qquad(k\ge0,\ 1\le r<p).
\end{equation}
For $1\le r<p$, divide $Rr$ by $p$ to write
\[
 Rr=pj(r)+\theta(r),\qquad
 j(r)=\left\lfloor\frac{Rr}{p}\right\rfloor,\qquad 1\le\theta(r)<p.
\]
Since $p\nmid R$, reduction modulo $p$ shows that $\theta$ is a permutation of $\{1,\ldots,p-1\}$. Set $K=kR+j(r)$ and $s=\beta_{[k]}(r)$. Then $R(pk+r)=pK+\theta(r)$ lies in the $p$-block of index $K$. By Step~5, its image $h(R(pk+r))=R(pk+s)=p(kR+j(s))+\theta(s)$ lies in the same $p$-block, so $j(s)=j(r)$.

For each $c\in\Z/L\Z$, define $\psi_c\colon\{1,\ldots,p-1\}\to\{1,\ldots,p-1\}$ by
\begin{equation}\label{eq:finite-reindex}
 \psi_c(\theta(r))=\theta(\beta_{[k]}(r)),
 \qquad\text{where }kR+j(r)\equiv c\pmod L.
\end{equation}
For each $r$, the residue class $[k]\in\Z/L\Z$ is uniquely determined because $\gcd(R,L)=1$, and since $\theta$ is a permutation, $\psi_c$ is well-defined on all of $\{1,\ldots,p-1\}$.

To see that $\psi_c$ is injective, take distinct $r_1,r_2\in\{1,\ldots,p-1\}$. For $i=1,2$, let $[k_i]$ be the unique solution of $kR+j(r_i)\equiv c\pmod L$ and set $s_i=\beta_{[k_i]}(r_i)$. If $[k_1]=[k_2]$, injectivity of $\beta_{[k_1]}$ gives $s_1\ne s_2$. If $[k_1]\ne[k_2]$, then $j(r_1)\ne j(r_2)$ (since equal values of $j$ would yield the same congruence), and since $j(s_i)=j(r_i)$ for $i=1,2$, we again have $s_1\ne s_2$. Injectivity of $\theta$ then yields $\psi_c(\theta(r_1))=\theta(s_1)\ne\theta(s_2)=\psi_c(\theta(r_2))$. Thus $\psi_c\in\Sym(\{1,\ldots,p-1\})$. Moreover, with $s=\beta_{[k]}(r)$ as above, Step~5 gives
\[
 \theta(s)-\theta(r)=R(s-r)
 =h(R(pk+r))-R(pk+r)\equiv0\pmod Y,
\]
so $\psi_c\in G_{p,Y}$.

Let $H$ be the self-homeomorphism defined by $(\psi_c)_{c\in\Z/L\Z}$ via~\eqref{eq:normal-map} (which belongs to $\Homeo(X_e)$ by Step~1). Since $R(pk+r)=p(kR+j(r))+\theta(r)$, equation~\eqref{eq:finite-reindex} gives $H(R(pk+r))=h(R(pk+r))$. On multiples of $Rp=P$, both $H$ and $h$ act as the identity. Hence $H$ and $h$ agree on $R\N$, so $H=h$ by $\mathsf B(S)$, establishing $\mathsf A(S)$ and completing the induction.

Finally, the values of $h$ on any $p$-block of index $K\in c$ determine $\sigma_c$ uniquely, so $h\mapsto(\sigma_c)_{c\in\Z/L\Z}$ is a bijection from $\Homeo(X_e)$ onto $\prod_{c\in\Z/L\Z}G_{p,Y}$. Since composition of maps of the form~\eqref{eq:normal-map} corresponds to componentwise composition, this bijection is a group isomorphism.
\end{proof}

The following elementary lemma determines when the fiber permutation group $G_{\ell,Y}$ is trivial and computes its order.

\begin{lemma}
\label{lem:fiber-permutation-threshold}
For a prime $\ell$ and $Y\ge1$,
\[
 G_{\ell,Y}=\{1\}
 \quad\Longleftrightarrow\quad
 Y\ge\ell-1.
\]
More generally, writing $m_c=|A_c^{\ell,Y}|$ for each $c\in\Z/Y\Z$, we have
\[
 |G_{\ell,Y}|=\prod_{c\in\Z/Y\Z}m_c!.
\]
\end{lemma}
\begin{proof}
If $Y\ge\ell-1$, distinct elements of $\{1,\ldots,\ell-1\}$ differ by at most $\ell-2$ and therefore lie in distinct residue classes modulo $Y$. If $Y<\ell-1$, then $1$ and $1+Y$ lie in the same fiber $A_{1\bmod Y}^{\ell,Y}$. Thus $G_{\ell,Y}$ is trivial if and only if $Y\ge\ell-1$. In general, each fiber $A_c^{\ell,Y}$ contributes $m_c!$ permutations independently.
\end{proof}

\begin{corollary}\label{cor:finite-rigidity}
For nonempty finite support $S$ with finite positive heights,
\[
 X_e\text{ is rigid}\quad\Longleftrightarrow\quad
 \prod_{q\in S\setminus\{\max S\}}q^{e(q)}\ge\max S-1.
\]
More explicitly, with $p=\max S$, $Y=\prod_{q\in S\setminus\{p\}}q^{e(q)}$, and $m_c=|\{1\le r<p:r\equiv c\pmod Y\}|$,
\[
 |\Homeo(X_e)|=
 \left(\prod_{c\in\Z/Y\Z}m_c!\right)^{p^{e(p)-1}}.
\]
\end{corollary}
\begin{proof}
Apply \cref{lem:fiber-permutation-threshold} to each of the $p^{e(p)-1}$ factors in~\eqref{eq:finite-homeo-group}.
\end{proof}

\begin{remark}\label{rem:empty-support}
When $S=\varnothing$, the only admissible common difference is $1$, so the specialization order is the usual order on $\N$. Since every order automorphism of $(\N,\le)$ is the identity, $X_e$ is rigid.
\end{remark}

\subsection{Exactly one infinite-height prime}\label{sec:one-infinite-height}

We now turn to the case $I_e\ne\varnothing$ (still assuming $S_e$ is finite). The first two lemmas apply for any nonempty $I_e$.

\begin{lemma}\label{lem:fm-displacement}
Let $h\in\Homeo(X_e)$ and define the displacement function $\delta\colon\N\to\Z$ by $\delta(n)=h(n)-n$. Fix $p\in I_e$. Then $h[p\N]=p\N$, and for every $x\notin p\N$ there exists $a\ge1$ such that $\delta(y)=\delta(x)$ for all $y\in C_e(x,p,a)$.
\end{lemma}
\begin{proof}
\noindent\textbf{Step 1.} We first show that divisibility by $p$ is preserved. Choose $N\ge1$ such that $p^N>2P_e$. Suppose for contradiction that $p\mid x$ while $p\nmid h(x)$. By continuity, there exists a basic neighborhood $x+d\Nzero$ with $h[x+d\Nzero]\subseteq h(x)+p^N\Nzero$. Since $p\mid x$ and $\gcd(x,d)=1$, we have $p\nmid d$, so $x+d\Nzero$ meets every residue class modulo $p^N$. On the other hand, for every $y\in x+d\Nzero$, the inclusion $h(y)\in h(x)+p^N\Nzero$ gives
\[
 y=h(y)-\delta(y)\equiv h(x)-\delta(y)\pmod{p^N}.
\]
By \cref{prop:common-finite-backbone}, $|\delta(y)|<P_e$, which allows at most $2P_e-1$ values of $\delta(y)$. Thus $x+d\Nzero$ can meet at most $2P_e-1$ residue classes modulo $p^N$, contradicting $p^N>2P_e$. Applying the same argument to $h^{-1}$ gives $h[p\N]=p\N$.

\smallskip\noindent\textbf{Step 2.} We next show that $\delta$ is locally constant on $p$-adic cylinders outside $p\N$. Fix $x\notin p\N$. By Step~1, $p\nmid h(x)$, so $h(x)+p^N\Nzero$ is a basic neighborhood of $h(x)$. By continuity and \cref{lem:tail-basis} (using $e(p)=\infty$ and $p\nmid x$ to ensure the exponent of $p$ is at least $N$), we may choose a basic neighborhood
\[
 U=x+p^a b\Nzero,\qquad h[U]\subseteq h(x)+p^N\Nzero,\qquad a\ge N,\quad\gcd(b,px)=1.
\]
Now fix any $y\in C_e(x,p,a)$, so $y\equiv x\pmod{p^a}$ and $p\nmid y$. By Step~1 and continuity, choose a basic neighborhood
\[
 V=y+p^m c\Nzero,\qquad h[V]\subseteq h(y)+p^N\Nzero,\qquad m\ge a,\quad\gcd(c,py)=1.
\]
Since $p\nmid bc$, the Chinese remainder theorem supplies $z\in\N$ satisfying
\[
 z\equiv y\pmod{p^m},\qquad z\equiv0\pmod r\quad\text{for every prime }r\mid bc.
\]
Because $m\ge a$, we also have $z\equiv x\pmod{p^a}$. For each prime $r\mid b$ we have $z\in r\N\subseteq r\N\cup(x+r^{v_r(b)}\Z)$, and similarly for each prime $r\mid c$. By \cref{prop:closure}, it follows that $z\in\overline U\cap\overline V$, while $z\equiv y\pmod p$ ensures $p\nmid z$.

By continuity and~\eqref{eq:closure},
\[
 h(z)\in\overline{h[U]}\subseteq\overline{h(x)+p^N\Nzero}=p\N\cup C_e(h(x),p,N).
\]
Since $h[p\N]=p\N$ and $p\nmid z$, we have $p\nmid h(z)$, forcing $h(z)\equiv h(x)\pmod{p^N}$. The same reasoning applied to $V$ gives $h(z)\equiv h(y)\pmod{p^N}$, so $h(y)\equiv h(x)\pmod{p^N}$. Since $y\equiv x\pmod{p^N}$ as well, we obtain
\[
 \delta(y)-\delta(x)=(h(y)-h(x))-(y-x)\equiv0\pmod{p^N}.
\]
Because $|\delta(y)-\delta(x)|<2P_e<p^N$, we conclude that $\delta(y)=\delta(x)$ for every $y\in C_e(x,p,a)$.
\end{proof}

\begin{lemma}\label{lem:fm-short-blocks}
Suppose $I_e\ne\varnothing$, and recall that $J_e=\prod_{p\in I_e}p$ and $\mathcal M_e=J_e\N$. Every $h\in\Homeo(X_e)$ fixes $\mathcal M_e$ pointwise and, for each $p\in I_e$, maps every interval between consecutive multiples of $p$,
\[
 D_k^{(p)}=\{pk+1,\ldots,pk+p-1\}\qquad(k\in\Nzero),
\]
onto itself.
\end{lemma}
\begin{proof}
Fix $h\in\Homeo(X_e)$ and put $\delta(n)=h(n)-n$.

\smallskip\noindent\textbf{Step 1.} Fix $p\in I_e$ and write $x=pk+a$ with $1\le a<p$. By \cref{lem:fm-displacement}, choose $N\ge1$ such that $\delta$ takes a constant value $c$ on $C=C_e(x,p,N)$. For each $0\le j<P_e/p$, since $\gcd(p^N,P_e)=p$ and $a+pj\equiv x\pmod p$, the Chinese remainder theorem yields $z_j\in C$ with $z_j\equiv a+pj\pmod{P_e}$. Write $z_j=t_jP_e+a+pj$ ($t_j\in\Nzero$). By \cref{prop:common-finite-backbone}, $h(z_j)=z_j+c\in B_{t_j}(P_e)$. Since $h$ fixes $P_e\N$ pointwise and $z_j\notin P_e\N$, injectivity rules out the right endpoint $(t_j+1)P_e$, so
\[
 1\le a+pj+c\le P_e-1.
\]
Taking $j=0$ and $j=P_e/p-1$ yields $1\le a+c\le p-1$, which means $h(x)\in D_k^{(p)}$. As this holds for every $x$ in the finite set $D_k^{(p)}$, injectivity gives $h(D_k^{(p)})=D_k^{(p)}$.

\smallskip\noindent\textbf{Step 2.} Let $x\in\mathcal M_e$, fix any $p\in I_e$, and set $d=Q_e(x)$. Because $p\mid x$ and $p\nmid d$, the least element of $A_x=(\mathord{\uparrow}_e x)\cap(\N\setminus p\N)$ is $x+d$, so the first interval $D_k^{(p)}$ met by $A_x$ has index
\[
 \min\{k:A_x\cap D_k^{(p)}\ne\varnothing\}
 =\frac{x}{p}+\left\lfloor\frac{d}{p}\right\rfloor.
\]
Since $h$ preserves $\preceq_e$ and $\N\setminus p\N$, we have $h[A_x]=A_{h(x)}$. By Step~1 and injectivity,
\[
 h[A_x\cap D_k^{(p)}]=h[A_x]\cap h[D_k^{(p)}]=A_{h(x)}\cap D_k^{(p)},
\]
so $A_x$ and $A_{h(x)}$ first meet $D_k^{(p)}$ at the same index $k$. Since $h(x)\in\mathcal M_e$ and $Q_e(h(x))=d$ by \cref{lem:fm-blocks}, we obtain
\[
 \frac{x}{p}+\left\lfloor\frac{d}{p}\right\rfloor
 =\frac{h(x)}{p}+\left\lfloor\frac{d}{p}\right\rfloor,
\]
whence $h(x)=x$.
\end{proof}

Now suppose $I_e=\{q\}$, and set $F=F_e=S_e\setminus\{q\}$, $Y=\prod_{r\in F}r^{e(r)}$, and $Q=Q_e$. Let $\Nzero^{(q)}$ denote $\Nzero$ equipped with the $q$-adic topology, in which the sets $\Nzero\cap(k+q^b\Z)$ ($b\in\Nzero$) form a neighborhood basis at $k$. Equip $G_{q,Y}$ with the discrete topology. A map $\Nzero^{(q)}\to G_{q,Y}$ is continuous if and only if it is locally constant; write $C(\Nzero^{(q)},G_{q,Y})$ for the group of all such maps under pointwise composition. To determine $\Homeo(X_e)$, we first obtain a family of block permutations from the invariance of the intervals $D_k^{(q)}$. Local constancy of the displacement and congruence modulo $Y$ show that this family belongs to $C(\Nzero^{(q)},G_{q,Y})$; we then verify the converse construction.

\begin{theorem}\label{thm:fm-one}
If $I_e=\{q\}$, then
\begin{equation}\label{eq:fm-group}
 \Homeo(X_e)\cong C(\Nzero^{(q)},G_{q,Y})
\end{equation}
as abstract groups. Under this isomorphism, a continuous map $\sigma\colon\Nzero^{(q)}\to G_{q,Y}$, $k\mapsto\sigma_k$, corresponds to the self-homeomorphism $h_\sigma\colon X_e\to X_e$ given by
\begin{equation}\label{eq:fm-normal}
 h_\sigma(qm)=qm\quad(m\ge1),\qquad
 h_\sigma(qk+a)=qk+\sigma_k(a)\quad(k\ge0,\ 1\le a<q).
\end{equation}
In particular, $X_e$ is rigid if and only if $Y\ge q-1$.
\end{theorem}
\begin{proof}
\noindent\textbf{Step 1.} Given $h\in\Homeo(X_e)$, \cref{lem:fm-short-blocks} shows that $h$ fixes $q\N$ pointwise and maps each $D_k^{(q)}$ onto itself. For each $k\in\Nzero$, define $\sigma_k\in\Sym(\{1,\ldots,q-1\})$ by $\sigma_k(a)=h(qk+a)-qk$; this is the unique permutation satisfying~\eqref{eq:fm-normal} on $D_k^{(q)}$.

\smallskip\noindent\textbf{Step 2.} We show that $\sigma_k\in G_{q,Y}$ for every $k$. Fix $x\notin q\N$. By \cref{lem:fm-displacement}, choose $t\ge1$ such that $\delta(n)=h(n)-n$ is constant on $C=C_e(x,q,t)$. Since $\gcd(q^t,Y)=1$, the Chinese remainder theorem gives $z\in C$ with $z\ge q$ and $z\equiv q\pmod Y$. Then $z\in\mathord{\uparrow}_e q=q+Y\Nzero$ (as $Q(q)=Y$). Since $h(q)=q$, this upper set is $h$-invariant, so $h(z)\equiv z\equiv q\pmod Y$. Therefore $Y\mid h(z)-z=\delta(x)$. Thus $h(n)\equiv n\pmod Y$ holds for all $n\notin q\N$, and trivially on $q\N$ as well:
\begin{equation}\label{eq:fm-congruence}
 h(n)\equiv n\pmod Y\qquad(n\in\N).
\end{equation}
In particular, $\sigma_k(a)-a=h(qk+a)-(qk+a)\equiv0\pmod Y$ for all $k\in\Nzero$ and $1\le a<q$, so $\sigma_k\in G_{q,Y}$.

\smallskip\noindent\textbf{Step 3.} We check that $k\mapsto\sigma_k$ is continuous on $\Nzero^{(q)}$. Fix $k\in\Nzero$ and $1\le a<q$. By \cref{lem:fm-displacement}, there exists $t\ge1$ such that $\delta$ is constant on $C_e(qk+a,q,t)$. Since $qj+a\in C_e(qk+a,q,t)$ if and only if $j\equiv k\pmod{q^{t-1}}$, the map $j\mapsto\sigma_j(a)-a=\delta(qj+a)$ is constant on a $q$-adic neighborhood of $k$. Intersecting these neighborhoods over $a=1,\ldots,q-1$ shows that $j\mapsto\sigma_j$ is locally constant, hence continuous.

\smallskip\noindent\textbf{Step 4.} Conversely, let $\sigma\colon\Nzero^{(q)}\to G_{q,Y}$ be continuous. Formula~\eqref{eq:fm-normal} defines a bijection $h_\sigma$ satisfying~\eqref{eq:fm-congruence}.

To prove continuity at $x=qk+a$ ($1\le a<q$), let $h_\sigma(x)+q^\alpha d\Nzero$ be a basic neighborhood of $h_\sigma(x)$, where $\alpha\ge0$, $d\mid Y$, and $\gcd(q^\alpha d,h_\sigma(x))=1$. Choose $L\ge0$ such that $\sigma_{k'}=\sigma_k$ whenever $k'\equiv k\pmod{q^L}$, and take $N\ge\max\{\alpha,L+1\}$. By~\eqref{eq:fm-congruence}, $\gcd(x,d)=1$, so $x+q^N d\Nzero$ is a basic neighborhood of $x$. Every $n\in x+q^N d\Nzero$ satisfies $n\ge x$ and, since $N\ge1$, has the form $n=qk'+a$ with $k'\ge k$ and $q^{N-1}\mid k'-k$. Hence $\sigma_{k'}=\sigma_k$, which gives
\[
 h_\sigma(n)-h_\sigma(x)=n-x\in q^N d\Nzero\subseteq q^\alpha d\Nzero.
\]

To prove continuity at $x=qm$, let $x+d\Nzero$ be a basic neighborhood of $h_\sigma(x)=x$, where $d\mid Y$ and $\gcd(x,d)=1$. For $n\in x+d\Nzero$,~\eqref{eq:fm-congruence} gives $h_\sigma(n)\equiv n\equiv x\pmod d$. Moreover, if $q\mid n$ then $h_\sigma(n)=n\ge x$; if $n=qk+a$ with $1\le a<q$, then $n\ge qm$ implies $k\ge m$, so $h_\sigma(n)=qk+\sigma_k(a)\ge x$. Thus $h_\sigma[x+d\Nzero]\subseteq x+d\Nzero$.

Since $G_{q,Y}$ is discrete, $k\mapsto\sigma_k^{-1}$ is also continuous, so the same argument shows that $h_\sigma^{-1}=h_{\sigma^{-1}}$ is continuous. Thus $h_\sigma\in\Homeo(X_e)$.

\smallskip\noindent\textbf{Step 5.} For any $\sigma,\tau\in C(\Nzero^{(q)},G_{q,Y})$, both $h_\sigma$ and $h_\tau$ fix $q\N$ pointwise, and $(h_\sigma\circ h_\tau)(qk+a)=qk+(\sigma_k\circ\tau_k)(a)$ for $1\le a<q$. Hence $\sigma\mapsto h_\sigma$ is a group isomorphism, proving~\eqref{eq:fm-group}. By \cref{lem:fiber-permutation-threshold}, $G_{q,Y}$ is trivial if and only if $Y\ge q-1$; when $G_{q,Y}$ is nontrivial, any constant nonidentity map $\sigma$ yields a nontrivial self-homeomorphism.
\end{proof}

\begin{remark}
If $S_e=\{q\}$, then $Y=1$, so $X_e$ is rigid if and only if $q=2$. Notice that the parameter maps in~\eqref{eq:fm-group} are locally constant functions on $\Nzero^{(q)}$ itself: they need neither extend continuously to the $q$-adic completion $\Z_q$ nor be globally periodic with respect to a single power of $q$.
\end{remark}

\subsection{At least two infinite-height primes}\label{sec:many-infinite-heights}

When at least two primes have infinite height, the Chinese remainder theorem couples their displacements and forces rigidity.

\begin{theorem}\label{thm:fm-many}
If $|I_e|\ge2$, then $X_e$ is rigid.
\end{theorem}
\begin{proof}
Write $I=I_e$, $J=\prod_{p\in I}p$, and $\mathcal M=J\N$. Let $h\in\Homeo(X_e)$ and put $\delta(n)=h(n)-n$.

We first show that $\delta$ is constant on $\N\setminus\mathcal M$. For any $p\in I$ and $x\notin p\N$, choose another prime $q\in I\setminus\{p\}$. By \cref{lem:fm-displacement}, there exist $a,b\ge1$ such that $\delta$ is constant on $C_e(x,p,a)$ and on $C_e(1,q,b)$, respectively. Since $p\ne q$, the Chinese remainder theorem produces $z\in\N$ with $z\equiv x\pmod{p^a}$ and $z\equiv1\pmod{q^b}$. Then $z$ lies in both cylinders, so $\delta(x)=\delta(z)=\delta(1)$. Because $\N\setminus\mathcal M=\bigcup_{p\in I}(\N\setminus p\N)$, we have $h(n)=n+\delta(1)$ for all $n\in\N\setminus\mathcal M$.

By \cref{lem:fm-short-blocks}, $h$ fixes $\mathcal M$ pointwise and therefore maps $\N\setminus\mathcal M$ onto itself. Since $1$ is the least element of $\N\setminus\mathcal M$, we obtain
\[
 1=\min h[\N\setminus\mathcal M]=1+\delta(1),
\]
so $\delta(1)=0$. Thus $h$ fixes $\N\setminus\mathcal M$ pointwise as well, and $h=\operatorname{id}_{\N}$.
\end{proof}

\section{Infinite support: coordinate reconstruction and cofinite rigidity}
\label{sec:cofinite-arbitrary}

Throughout this section, write $S=S_e$, $T=\PP\setminus S$, $F=F_e$, and $I=I_e$, and assume that $S$ is infinite. The additional assumption that $T$ is finite will be imposed only after \cref{lem:external-finite-rank} when we prove rigidity.

We first show that for infinite support, the unit residue classes generate the full topology $\taue$ without needing upward rays.

\begin{proposition}
\label{prop:inf-terminal-basis}
The following cylinders, with $a\in\N$, form a subbasis for $\taue$:
\begin{equation}\label{eq:inf-terminal-generators}
 \begin{aligned}
 &C_e(a,p,k),&&p\in I,\ k\ge1,\ p\nmid a,\\
 &C_e(a,q,e(q)),&&q\in F,\ q\nmid a.
 \end{aligned}
\end{equation}
\end{proposition}
\begin{proof}
\noindent\textbf{Step 1.} Let $\rho_e$ be the topology generated by~\eqref{eq:inf-terminal-generators}. We first show that every finite-height cylinder is $\rho_e$-open. For $q\in F$, put $m=e(q)$. If $1\le k<m$ and $q\nmid a$, then
\[
 C_e(a,q,k)=\bigsqcup_{t=0}^{q^{m-k}-1} C_e(a+tq^k,q,m),
\]
because $q^m\mid n-(a+tq^k)$ for a unique $0\le t<q^{m-k}$ whenever $q^k\mid n-a$, and each $a+tq^k$ is coprime to $q$. All infinite-height cylinders already appear in~\eqref{eq:inf-terminal-generators}.

\smallskip\noindent\textbf{Step 2.} Each cylinder $C_e(a,p,k)$ in~\eqref{eq:inf-terminal-generators} is the union of the admissible progressions $n+p^k\Nzero$ over $n\in C_e(a,p,k)$, so $\rho_e\subseteq\taue$. Conversely, for any distinct $x,y\in\N$, choosing $p\in S$ with $p>\max\{x,y\}$ gives disjoint $\rho_e$-open cylinders $C_e(x,p,1)$ and $C_e(y,p,1)$ separating $x$ and $y$. Thus $\rho_e$ is Hausdorff, and all finite subsets of $\N$ are $\rho_e$-closed. Removing the finitely many elements below $a$ from $C_e(a,p,k)$ shows that $a+p^k\Nzero$ is $\rho_e$-open for all $k\ge1$, and the upward rays $a+\Nzero$ ($k=0$) are $\rho_e$-open because their complements are finite. Hence $\taue\subseteq\rho_e$, so $\rho_e=\taue$.
\end{proof}

\begin{remark}\label{rem:inf-golomb-cylinders}
Golomb used this description by truncations of coprime progressions to show that the set of $m\in\N$ for which $6m-1$ and $6m+1$ are both prime is closed~\cite[p.~181]{Golomb1962}. For the Golomb and Kirch cases, see also~\cite[Remark~2.18]{Longhi2023}, which adapts~\cite[Proposition~2.1]{Orum2016}.
\end{remark}

Following~\cite[\S2 and Proposition~1]{Banakh2021b}, we now reconstruct the action of a self-homeomorphism on active primes using filters of neighborhood closures.

When $e(2)=1$, the only unit residue class modulo $2$ is the odd numbers, so its closure is all of $\N$. We set $2$ aside temporarily by defining
\begin{equation}\label{eq:S-sharp}
 S^\sharp=
 \begin{cases}
 S\setminus\{2\},&2\in S\text{ and }e(2)=1,\\
 S,&\text{otherwise}.
 \end{cases}
\end{equation}
We will recover the invariance of the parity classes separately in \cref{lem:cof-parity-repair}.

For a nonempty finite set $E\subseteq\N$, let $\mathscr F_E$ be the family of all $B\subseteq\N$ for which there exist neighborhoods $V_x\ni x$ ($x\in E$) such that
\begin{equation}\label{eq:anchored-filter-global}
 \bigcap_{x\in E}\overline{V_x}\subseteq B.
\end{equation}
Similarly, put
\begin{equation}\label{eq:unanchored-filter-global}
 \mathscr F_*=\left\{B\subseteq\N:
 \bigcap_{i=1}^m\overline{V_i}\subseteq B
 \text{ for some nonempty open }V_1,\ldots,V_m\right\}.
\end{equation}
Both families are upward closed and omit $\varnothing$ by \cref{cor:superconnected}. Intersecting the corresponding neighborhoods at each $x\in E$ shows that $\mathscr F_E$ is closed under finite intersections, and combining finite witness families shows the same for $\mathscr F_*$. Thus both are filters. Every $h\in\Homeo(X_e)$ transports $\mathscr F_E$ to $\mathscr F_{h[E]}$ and preserves $\mathscr F_*$.

For a nonempty finite set $E\subseteq\N$ and $p\in S^\sharp$, define the \emph{local state} $s_p(E)$ of $E$ at $p$ as follows:
\begin{enumerate}
\item Set $s_p(E)=\mathsf c$ if $p$ divides every element of $E$.
\item If $p\in F$ and the elements of $E$ not divisible by $p$ are nonempty and all congruent modulo $p^{e(p)}$, set $s_p(E)=u$, where $u\in(\Z/p^{e(p)}\Z)^\times$ is their common residue class.
\item If $p\in I$ and $x$ is the unique element of $E$ not divisible by $p$, set $s_p(E)=x$.
\item Set $s_p(E)=\mathsf a$ otherwise (that is, if $p\in F$ and two elements of $E$ coprime to $p$ are incongruent modulo $p^{e(p)}$, or if $p\in I$ and at least two elements of $E$ are coprime to $p$).
\end{enumerate}

At an odd prime $p$ of height $1$, this recovers the three Kirch states from~\cite[\S3, Lemmas~5--6]{Banakh2021b}: all elements divisible by $p$ ($\mathsf c$), all $p$-coprime elements in one residue class ($u$), or at least two residue classes ($\mathsf a$). At infinite height, any two distinct $p$-units $x,y\in E$ are separated modulo $p^k$ once $p^k>|x-y|$, so an intermediate state occurs only when $E$ has a unique $p$-unit $x$. Note also that if $|E|\ge2$ and $p>\max E$, any two distinct elements of $E$ are $p$-units incongruent modulo $p$, so $s_p(E)=\mathsf a$.

\begin{lemma}
\label{lem:cof-local-state}
For $|E|\ge2$, a filter basis for $\mathscr F_E$ consists of all finite intersections formed by choosing, at a finite subset of $S^\sharp$:
\begin{enumerate}
\item $p\N$ at primes with $s_p(E)=\mathsf a$;
\item $\overline{u+p^{e(p)}\Nzero}$ at primes with finite-height state $s_p(E)=u$;
\item $\overline{x+p^k\Nzero}$ (for any $k\ge1$) at primes with infinite-height state $s_p(E)=x$;
\end{enumerate}
and choosing no set at primes with state $\mathsf c$. Order the local states by $\mathsf c<u<\mathsf a$ at finite height and $\mathsf c<x<\mathsf a$ at infinite height, with distinct intermediate states incomparable. Then, for $|E|,|D|\ge2$,
\[
 \mathscr F_E\subseteq\mathscr F_D
 \quad\Longleftrightarrow\quad
 s_p(E)\le s_p(D)\text{ for all }p\in S^\sharp.
\]
\end{lemma}
\begin{proof}
\noindent\textbf{Step 1.} We verify that the listed intersections form a basis for $\mathscr F_E$. Given such an intersection over a finite set $J\subseteq S^\sharp$, define an exponent $k_p\le e(p)$ for each $p\in J$ as follows: use the exponent $k_p$ of the chosen closure if $s_p(E)$ is intermediate; if $s_p(E)=\mathsf a$, take $k_p=e(p)$ when $p\in F$, and choose $k_p\ge1$ with $p^{k_p}>|x-y|$ for two $p$-units $x,y\in E$ when $p\in I$. For each $x\in E$, put
\[
 d_x=\prod_{\substack{p\in J\\p\nmid x}}p^{k_p}.
\]
Then $d_x\in\De$ and $\gcd(x,d_x)=1$. By \cref{prop:closure}, at an intermediate state the closures $\overline{x+p^{k_p}\Nzero}$ for those $x\in E$ with $p\nmid x$ all agree (while elements divisible by $p$ impose no condition at $p$), whereas at state $\mathsf a$ two of the unit classes are disjoint and every closure contains $p\N$. Thus the listed intersection equals $\bigcap_{x\in E}\overline{x+d_x\Nzero}\in\mathscr F_E$.

Conversely, take arbitrary basic neighborhoods $x+b_x\Nzero$ ($x\in E$). Let $J\subseteq S^\sharp$ consist of all primes with intermediate state for $E$ together with all primes in $S^\sharp$ dividing some $b_x$. For each $p\in J$, take $k_p=e(p)$ if $p\in F$, and if $p\in I$ choose a finite exponent $k_p\ge\max\{1,\max_{x\in E}v_p(b_x)\}$ large enough to satisfy the separation condition above whenever $s_p(E)=\mathsf a$. Since each $b_x\in\De$, we have $v_p(b_x)\le e(p)$ for $p\in F$ and $e(p)=\infty$ for $p\in I$, so this choice is admissible, and for each $p\in S^\sharp$ dividing $b_x$ we have $p\nmid x$ and $p^{v_p(b_x)}\mid p^{k_p}$, so $\overline{x+p^{k_p}\Nzero}\subseteq\overline{x+p^{v_p(b_x)}\Nzero}$. If $2\notin S^\sharp$, the only possible closure at $2$ is $\overline{x+2\Nzero}=\N$ for odd $x$. Hence the resulting basis intersection is contained in $\bigcap_{x\in E}\overline{x+b_x\Nzero}$.

\smallskip\noindent\textbf{Step 2.} Assume $s_p(E)\le s_p(D)$ for all $p\in S^\sharp$. Every factor at $p$ in a basic member of $\mathscr F_E$ either agrees with a basic factor of $\mathscr F_D$ (when $s_p(E)=s_p(D)$ is intermediate) or contains $p\N\in\mathscr F_D$ (when $s_p(D)=\mathsf a$). Hence $\mathscr F_E\subseteq\mathscr F_D$.

\smallskip\noindent\textbf{Step 3.} Conversely, suppose $s_p(E)\nleq s_p(D)$ for some $p\in S^\sharp$. We construct $H\in\mathscr F_E$ such that $B\setminus H\ne\varnothing$ for every basic member $B\in\mathscr F_D$.
\begin{enumerate}[label=(\alph*)]
\item If $s_p(E)=\mathsf a$, set $H=p\N$. If $B$ uses $\overline{v+p^j\Nzero}$ at $p$, require $n\equiv v\pmod{p^j}$; if $s_p(D)=\mathsf c$, require $n\equiv1\pmod p$.
\item If $s_p(E)$ and $s_p(D)$ are distinct intermediate states represented by $u$ and $v$, choose $k\le e(p)$ with $u\not\equiv v\pmod{p^k}$ (taking $k=e(p)$ if $p\in F$) and set $H=\overline{u+p^k\Nzero}$. If $B$ uses $\overline{v+p^j\Nzero}$ at $p$, require $n\equiv v\pmod{p^{\max\{j,k\}}}$.
\item If $s_p(E)$ is intermediate, represented by $u$, and $s_p(D)=\mathsf c$, set $H=\overline{u+p^k\Nzero}$ with $k=e(p)$ for $p\in F$ and $k=2$ for $p\in I$. Since $p\in S^\sharp$ at finite height and $\varphi(p^2)\ge2$ at infinite height, there exists a unit class $w\not\equiv u\pmod{p^k}$; require $n\equiv w\pmod{p^k}$.
\end{enumerate}
If $B$ has no factor at $p$, impose $n\equiv1\pmod p$ in case (a), and $n\equiv v\pmod{p^k}$ in case (b). In each case, the congruence at $p$ satisfies the condition of $B$ at $p$ (if any) while ensuring $n\notin H$ by~\eqref{eq:closure}. For every other prime $q$ appearing in $B$, require $n\equiv0\pmod q$ (as $q\N$ is contained in the factor of $B$ at $q$). By the Chinese remainder theorem, there exists $n\in B\setminus H$, so $H\notin\mathscr F_D$.
\end{proof}

For a singleton $\{z\}$, the same calculation gives a filter basis for $\mathscr F_{\{z\}}$ consisting of finite intersections of $\overline{z+p^{e(p)}\Nzero}$ ($p\in F\cap S^\sharp$, $p\nmid z$) and $\overline{z+p^k\Nzero}$ ($p\in I$, $p\nmid z$, $k\ge1$).

Consider the partially ordered set
\[
 \mathscr P_e=\{\mathscr F_E:2\le|E|<\infty\}\cup\{\mathscr F_*\},
\]
ordered by inclusion. A \emph{coatom} of $\mathscr P_e$ is an element maximal among those properly contained in $\mathscr F_*$ (analogous to~\cite[Lemmas~11--12]{Banakh2021b}). For $p\in F\cap S^\sharp$ and $u\in(\Z/p^{e(p)}\Z)^\times$, let $\mathscr M(p,u)$ be the filter generated by $\overline{u+p^{e(p)}\Nzero}$ together with $q\N$ for all $q\in S^\sharp\setminus\{p\}$. For $p\in I$ and $x\in\N$ with $p\nmid x$, let $\mathscr M(p,x)$ be the filter generated by $\overline{x+p^k\Nzero}$ ($k\ge1$) together with $q\N$ for all $q\in S^\sharp\setminus\{p\}$.

\begin{lemma}
\label{lem:cof-coatoms-incidence}
The coatoms of $\mathscr P_e$ are precisely
\begin{equation}\label{eq:cof-coatom-list}
 \begin{aligned}
 &\mathscr M(p,u),&&p\in F\cap S^\sharp,\quad u\in(\Z/p^{e(p)}\Z)^\times,\\
 &\mathscr M(p,x),&&p\in I,\quad x\in\N,\quad p\nmid x,
 \end{aligned}
\end{equation}
and each is realized as $\mathscr F_E$ for some finite $E$ with $|E|\ge2$. Moreover, the incidence map $A(\mathscr M)=\{z\in\N:\mathscr F_{\{z\}}\subseteq\mathscr M\}$ satisfies
\begin{equation}\label{eq:cof-incidence-formula}
 \begin{aligned}
 A(\mathscr M(p,u))&=p\N\cup C_e(u,p,e(p)),&&p\in F\cap S^\sharp,\\
 A(\mathscr M(p,x))&=p\N\cup\{x\},&&p\in I.
 \end{aligned}
\end{equation}
\end{lemma}
\begin{proof}
\noindent\textbf{Step 1.} We first note that $\mathscr F_*$ has basis $\{\bigcap_{p\in A}p\N:A\subseteq S^\sharp\text{ finite}\}$, corresponding to state $\mathsf a$ at every $p\in S^\sharp$. Indeed, for each $p\in S^\sharp$ we can choose an admissible exponent $k\ge1$ and $p$-units $u\not\equiv v\pmod{p^k}$, so $\overline{u+p^k\Nzero}\cap\overline{v+p^k\Nzero}=p\N\in\mathscr F_*$. Conversely, every finite intersection of closures of nonempty open sets contains $\bigcap_{p\in A}p\N$, where $A\subseteq S^\sharp$ is the finite set of primes dividing the common differences (an omitted prime $2$ contributes only $\overline{u+2\Nzero}=\N$).

\smallskip\noindent\textbf{Step 2.} To realize each filter in~\eqref{eq:cof-coatom-list}, take a positive representative $y$ of $u$ for $\mathscr M(p,u)$, or $y=x$ for $\mathscr M(p,x)$, and set $E=\{y,p,2p,3p\}$. Then $E$ has the required intermediate state at $p$. For every odd $q\in S^\sharp\setminus\{p\}$, the elements $p,2p\in E$ are coprime to $q$ and distinct modulo $q$; if $2\in S^\sharp$ and $p\ne2$, then $e(2)\ge2$ and $p,3p$ are odd and distinct modulo $4$. Thus $s_q(E)=\mathsf a$ for every $q\in S^\sharp\setminus\{p\}$, so $\mathscr F_E$ equals the given filter.

\smallskip\noindent\textbf{Step 3.} For $D_0=\{1,2,3\}$, the state is $\mathsf a$ at every prime of $S^\sharp$: use $1,2$ at odd primes and $1,3$ modulo $4$ at a retained prime $2$. Hence $\mathscr F_{D_0}=\mathscr F_*$ by Step~1 and \cref{lem:cof-local-state}, so the local-state comparison also applies to the top filter. By \cref{lem:cof-local-state}, any $\mathscr F_E\subsetneq\mathscr F_*$ has $s_p(E)<\mathsf a$ at some $p\in S^\sharp$. Raising its state at $p$ to an intermediate state if necessary and setting the state at all other primes in $S^\sharp$ to $\mathsf a$ embeds $\mathscr F_E$ into one of the filters in~\eqref{eq:cof-coatom-list}. Conversely, no element of $\mathscr P_e$ lies strictly between such a filter and $\mathscr F_*$, since only the state at $p$ can increase. Thus~\eqref{eq:cof-coatom-list} lists all coatoms of $\mathscr P_e$.

\smallskip\noindent\textbf{Step 4.} Fix $z\in\N$ and a coatom $\mathscr M$ at $p\in S^\sharp$, with $w$ denoting either a positive representative of $u$ (if $\mathscr M=\mathscr M(p,u)$) or $x$ (if $\mathscr M=\mathscr M(p,x)$). For any $q\in S^\sharp\setminus\{p\}$, every basis generator of $\mathscr F_{\{z\}}$ at $q$ contains $q\N\in\mathscr M$. Hence, if $p\mid z$, then $\mathscr F_{\{z\}}\subseteq\mathscr M$.

Assume $p\nmid z$. For each level $k$ appearing in the basis of $\mathscr F_{\{z\}}$ at $p$ ($k=e(p)$ if $p\in F$, or any $k\ge1$ if $p\in I$), we claim that
\[
 \overline{z+p^k\Nzero}\in\mathscr M
 \quad\Longleftrightarrow\quad z\equiv w\pmod{p^k}.
\]
The backward implication is immediate since $\overline{z+p^k\Nzero}=\overline{w+p^k\Nzero}\in\mathscr M$. If $z\not\equiv w\pmod{p^k}$, any basic member of $\mathscr M$ can be refined to include a closure $\overline{w+p^j\Nzero}$ at $p$ (as these closures are nested in $j$). By the Chinese remainder theorem, there exists $n\in\N$ such that $n\equiv w\pmod{p^{\max\{j,k\}}}$ and $n\equiv0\pmod q$ for every other prime $q$ used in that basic member. Then $n$ lies in the basic member, yet $p\nmid n$ and $n\not\equiv z\pmod{p^k}$, so $n\notin\overline{z+p^k\Nzero}$. Thus $\overline{z+p^k\Nzero}\notin\mathscr M$.

Consequently, for $p\nmid z$, we have $\mathscr F_{\{z\}}\subseteq\mathscr M$ if and only if $z\equiv u\pmod{p^{e(p)}}$ (when $p\in F$) or $z\equiv x\pmod{p^k}$ for all $k\ge1$, i.e., $z=x$ (when $p\in I$). Combining this with $p\mid z$ yields~\eqref{eq:cof-incidence-formula}.
\end{proof}

For each $p\in S^\sharp$, group the coatoms associated with $p$ into the family
\[
 \mathscr B_p=
 \begin{cases}
  \{\mathscr M(p,u):u\in(\Z/p^{e(p)}\Z)^\times\},&p\in F\cap S^\sharp,\\
  \{\mathscr M(p,x):x\in\N,\ p\nmid x\},&p\in I.
 \end{cases}
\]
The next lemma characterizes the families $\mathscr B_p$ purely in terms of the poset $\mathscr P_e$ and the incidence map $A$.

\begin{lemma}
\label{lem:cof-sunflower}
Let $\mathscr B$ be a set of coatoms of $\mathscr P_e$ with $|\mathscr B|\ge2$. Then $\mathscr B=\mathscr B_p$ for some $p\in S^\sharp$ if and only if:
\begin{enumerate}
\item there is a set $D_{\mathscr B}\subseteq\N$ such that $A(\mathscr M)\cap A(\mathscr M')=D_{\mathscr B}$ for all distinct $\mathscr M,\mathscr M'\in\mathscr B$; and
\item $\bigcup_{\mathscr M\in\mathscr B}A(\mathscr M)=\N$.
\end{enumerate}
In that case, $D_{\mathscr B}=p\N$.
\end{lemma}
\begin{proof}
For fixed $p\in S^\sharp$, the sets $A(\mathscr M)\setminus p\N$ ($\mathscr M\in\mathscr B_p$) partition $\N\setminus p\N$ (into unit residue classes modulo $p^{e(p)}$ if $p\in F$, or into singletons if $p\in I$). Hence the sets $A(\mathscr M)$ have pairwise intersection $p\N$ and union $\N$, proving necessity.

Conversely, suppose (i) and (ii) hold. Assume for contradiction that $\mathscr B$ contains members of $\mathscr B_p$ and $\mathscr B_q$ for distinct primes $p,q$, with incidence sets $p\N\cup L$ and $q\N\cup K$, respectively. There cannot be a member of $\mathscr B_r$ for any $r\notin\{p,q\}$: its incidence set would intersect $p\N\cup L$ in a set containing $pr\N$, whereas the Chinese remainder theorem supplies a multiple of $pr$ outside $q\N\cup K$, contradicting (i). Likewise, $\mathscr B$ cannot contain a second member of $\mathscr B_p$: two members of $\mathscr B_p$ intersect in $p\N$, which is not contained in $q\N\cup K$. The same holds for $\mathscr B_q$, so $|\mathscr B|=2$. However, by the Chinese remainder theorem, there is an integer coprime to both $p$ and $q$ lying outside $L\cup K$ (choosing it sufficiently large if $L$ or $K$ is a singleton), contradicting (ii).

Thus $\mathscr B\subseteq\mathscr B_p$ for a single prime $p\in S^\sharp$, and (ii) then forces $\mathscr B=\mathscr B_p$.
\end{proof}

\begin{theorem}
\label{thm:active-coordinate-reconstruction}
Every $h\in\Homeo(X_e)$ induces a permutation $\sigma\colon S^\sharp\to S^\sharp$ such that
\begin{equation}\label{eq:cof-multiple-transport}
 h[p\N]=\sigma(p)\N\qquad(p\in S^\sharp).
\end{equation}
For each $p\in S^\sharp$, either $p,\sigma(p)\in F$ or $p,\sigma(p)\in I$. Moreover, for $p\in F\cap S^\sharp$ and $q=\sigma(p)$, there is a bijection $\theta_p\colon(\Z/p^{e(p)}\Z)^\times\to(\Z/q^{e(q)}\Z)^\times$ such that
\begin{equation}\label{eq:cof-terminal-transport}
 h[C_e(u,p,e(p))]=C_e(\theta_p(u),q,e(q)),
\end{equation}
and in particular
\begin{equation}\label{eq:cof-totient-cardinality}
 \varphi(p^{e(p)})=\varphi(q^{e(q)}).
\end{equation}
\end{theorem}
\begin{proof}
For any filter $\mathscr F$ on $\N$, set $h_*\mathscr F=\{h[B]:B\in\mathscr F\}$. Because $h_*\mathscr F_E=\mathscr F_{h[E]}$ and $h_*\mathscr F_*=\mathscr F_*$, $h_*$ is an order automorphism of $\mathscr P_e$. Since $h[A(\mathscr M)]=A(h_*\mathscr M)$, $h_*$ preserves conditions~(i) and~(ii) of \cref{lem:cof-sunflower}. Hence for each $p\in S^\sharp$ there is a unique $\sigma(p)\in S^\sharp$ such that $h_*[\mathscr B_p]=\mathscr B_{\sigma(p)}$. Applying the same argument to $h^{-1}$ shows that $\sigma$ is a permutation of $S^\sharp$. Comparing common intersections $D_{\mathscr B_p}$ gives~\eqref{eq:cof-multiple-transport}.

Since $\mathscr B_p$ has $\varphi(p^{e(p)})$ elements when $p\in F$ and infinitely many when $p\in I$, $\sigma$ preserves $F\cap S^\sharp$ and $I$. For $p\in F\cap S^\sharp$ and $q=\sigma(p)$, define $\theta_p(u)$ by $h_*\mathscr M(p,u)=\mathscr M(q,\theta_p(u))$. Then $\theta_p$ is bijective, and removing $p\N$ and $q\N$ from $h[A(\mathscr M(p,u))]=A(\mathscr M(q,\theta_p(u)))$ yields~\eqref{eq:cof-terminal-transport}.
\end{proof}

For $p\in I$, we next adapt the characterization of the $p$-adic topology on $U_p=\N\setminus p\N$ via clopen sets~\cite[Lemmas~6--7]{Banakh2021a} (see also~\cite[Proposition~3.1 and Theorem~3.3]{Spirito2020}).

\begin{proposition}
\label{prop:common-padic-reflection}
For $p\in I$, the topology generated by the clopen subsets of $U_p\subseteq X_e$ is its ordinary $p$-adic topology. Consequently, if $h\in\Homeo(X_e)$ satisfies $h[p\N]=q\N$ with $p,q\in I$, then $h|_{U_p}\colon U_p\to U_q$ is a homeomorphism when $U_p$ and $U_q$ carry the $p$-adic and $q$-adic topologies, respectively.
\end{proposition}
\begin{proof}
Every unit residue class modulo $p^k$ ($k\ge1$) is clopen in $U_p$, since both the class and its complement in $U_p$ are unions of admissible progressions.

Conversely, let $\Omega\subseteq U_p$ be clopen and take $x\in\Omega$. Choose an admissible progression $x+p^n b\Nzero\subseteq\Omega$ with $n\ge1$ and $p\nmid b$. If some $y\in U_p\setminus\Omega$ satisfied $y\equiv x\pmod{p^n}$, we could choose an admissible progression $y+p^m c\Nzero\subseteq U_p\setminus\Omega$ with $m\ge n$ and $p\nmid c$. By the Chinese remainder theorem, there exists $z\in\N$ such that $z\equiv y\pmod{p^m}$ (hence $z\equiv x\pmod{p^n}$) and $z\equiv0\pmod r$ for every prime $r\mid bc$. By \cref{prop:closure},
\[
 z\in\overline{x+p^n b\Nzero}\cap\overline{y+p^m c\Nzero}.
\]
Since $p\nmid z$, intersecting with $U_p$ places $z$ in both relatively closed sets $\Omega$ and $U_p\setminus\Omega$, a contradiction. Thus $U_p\cap(x+p^n\Z)\subseteq\Omega$. Because homeomorphisms preserve relatively clopen sets, the final assertion follows.
\end{proof}

Fix $h\in\Homeo(X_e)$ and its induced permutation $\sigma$. Combining \cref{prop:common-padic-reflection} with \cref{thm:active-coordinate-reconstruction} gives the following corollary.

\begin{corollary}
\label{cor:cof-padic-reflection}
For $p\in I$ and $q=\sigma(p)$, the restriction $h|_{U_p}\colon U_p\to U_q$ is a homeomorphism from the $p$-adic topology on $U_p$ to the $q$-adic topology on $U_q$.
\end{corollary}

We now restore the prime $2$ when $e(2)=1$.

\begin{lemma}
\label{lem:cof-parity-repair}
If $2\in S$ and $e(2)=1$, then $h[2\N]=2\N$ and $h[1+2\Nzero]=1+2\Nzero$. Thus $\sigma$ extends to a permutation of $S$ by setting $\sigma(2)=2$.
\end{lemma}
\begin{proof}
Let $\rho_o$ be the topology on $\N$ generated by all admissible unit residue classes at odd active primes. By \cref{prop:inf-terminal-basis} (Step~1), each finite-height class is a union of highest-level classes, which $h$ and $h^{-1}$ transport by \cref{thm:active-coordinate-reconstruction}; at infinite height, \cref{cor:cof-padic-reflection} gives continuity in both directions. Thus $h$ is a $\rho_o$-homeomorphism. Every nonempty basic $\rho_o$-open set contains an arithmetic progression of odd common difference, so both $O=1+2\Nzero$ and $E=2\N$ are $\rho_o$-dense.

Suppose for contradiction that the $\taue$-open set $h[O]$ contains an even integer $x$. Then $h[O]$ contains an admissible progression $x+d\Nzero$, where $d$ is odd since $\gcd(x,d)=1$. Defining $e_o(2)=0$ and $e_o(p)=e(p)$ for odd $p$, we see from \cref{prop:inf-terminal-basis} that $\tau_{e_o}=\rho_o$. Because $x+d\Nzero$ is admissible for $e_o$, it is $\rho_o$-open, contradicting the fact that $h[E]$ is $\rho_o$-dense and disjoint from $h[O]$. Hence $h[O]\subseteq O$. Applying the same argument to $h^{-1}$ gives $h[O]=O$, and taking complements gives $h[E]=E$.
\end{proof}

We have thus obtained, for every infinite support $S$, a permutation $\sigma\colon S\to S$ preserving $F$ and $I$ and satisfying $h[p\N]=\sigma(p)\N$ for all $p\in S$. Setting $\Pi_S(n)=\{p\in S:p\mid n\}$, this reads
\begin{equation}\label{eq:active-support-transport}
 \Pi_S(h(n))=\sigma[\Pi_S(n)]\qquad(n\in\N).
\end{equation}
(For the Golomb case, see~\cite[Theorem~5.1(3)]{Banakh2018}.) For any set $A\subseteq\PP$, define
\begin{equation}\label{eq:T-units}
 \Gamma_A=\langle a:a\in A\rangle\le\Q_{>0}^{\times},\qquad
 \mathcal U_T=\Gamma_T\cap\N,\qquad \mathcal Y_T=\N\setminus\mathcal U_T.
\end{equation}
Thus $\mathcal U_T$ consists of the positive integers whose prime divisors all lie in $T$, while $\mathcal Y_T$ consists of those having at least one prime divisor in $S$. Equation~\eqref{eq:active-support-transport} implies
\begin{equation}\label{eq:T-unit-invariance}
 h[\mathcal U_T]=\mathcal U_T,\qquad h[\mathcal Y_T]=\mathcal Y_T.
\end{equation}
At this stage we only know that $h(1)\in\mathcal U_T$.

\begin{lemma}
\label{lem:connected-active-layer}
The subspace $\mathcal Y_T\subseteq X_e$ is dense and superconnected, and hence connected.
\end{lemma}
\begin{proof}
Every basic progression $a+d\Nzero$ contains a point divisible by some $q\in S$ with $q\nmid d$ (by the Chinese remainder theorem), so $\mathcal Y_T$ is dense in $X_e$.

To prove superconnectedness, write $\rad(d)=\prod_{p\mid d}p$. We first show that
\begin{equation}\label{eq:Y-backbone}
 \rad(d)\N\cap \mathcal Y_T\subseteq\overline{(a+d\Nzero)\cap \mathcal Y_T}^{\,\mathcal Y_T}
\end{equation}
for every basic progression $a+d\Nzero$. Take $z\in\rad(d)\N\cap \mathcal Y_T$ and a relative basic neighborhood $(z+c\Nzero)\cap \mathcal Y_T$. Since every prime divisor of $d$ divides $z$ while $\gcd(z,c)=1$, we have $\gcd(c,d)=1$. Choosing $q\in S$ with $q\nmid cd$, the congruences $w\equiv z\pmod c$, $w\equiv a\pmod d$, and $w\equiv0\pmod q$ have arbitrarily large positive solutions $w$, which lie in $(z+c\Nzero)\cap(a+d\Nzero)\cap \mathcal Y_T$. This proves~\eqref{eq:Y-backbone}.

Given finitely many nonempty relatively open subsets of $\mathcal Y_T$, choose relative basic progressions $(a_i+d_i\Nzero)\cap \mathcal Y_T$ inside them, and let $D$ be the product of the primes dividing any $d_i$. By~\eqref{eq:Y-backbone}, all their relative closures contain the nonempty set $D\N\cap \mathcal Y_T$. Thus $\mathcal Y_T$ is superconnected, and hence connected.
\end{proof}

Next, we record two local-global principles for finitely generated subgroups of $\Q_{>0}^{\times}$, derived from Perucca's results~\cite[Theorem~8 and Lemma~10]{Perucca2010} (see also~\cite[Theorem~1]{CorralesRodriganez1997} for cyclic subgroups).

For a positive rational $u=a/b$ in lowest terms and a prime $r\nmid ab$, denote its reduction modulo $r$ by $\overline u_r=[a]_r[b]_r^{-1}\in\mathbb F_r^\times$, where $c\mapsto[c]_r$ is the residue map $\Z\to\mathbb F_r$. For a finitely generated subgroup $\Gamma\le\Q_{>0}^{\times}$, reduction defines a homomorphism $\Gamma\to\mathbb F_r^\times$ for all primes $r$ dividing neither the numerators nor the denominators of a finite generating set, with image $\overline\Gamma_r=\{\overline u_r:u\in\Gamma\}\le\mathbb F_r^\times$.

\begin{lemma}\label{lem:external-finite-rank}
Let $\Gamma\le\mathbb Q_{>0}^{\times}$ be finitely generated.
\begin{enumerate}
\item Let $c\in\mathbb Q_{>0}^{\times}$. If $\overline c_r\in\overline\Gamma_r$ for all but finitely many primes $r$ for which both $\overline c_r$ and $\overline\Gamma_r$ are defined, then $c\in\Gamma$.
\item Let $b\in\mathbb Q_{>0}^{\times}$, and suppose that $b^n\notin\Gamma$ for every $n\in\N$. Then, for every prime $\ell$ and every finite set $E$ of primes, there exists a prime $r\notin E$ for which both $\overline b_r$ and $\overline\Gamma_r$ are defined and such that the coset $\overline b_r\,\overline\Gamma_r\in\mathbb F_r^\times/\overline\Gamma_r$ has order $d_r$ divisible by $\ell$ (equivalently, $\ell\mid d_r$, where $d_r\ge1$ is minimal with $\overline b_r^{\,d_r}\in\overline\Gamma_r$).
\end{enumerate}
\end{lemma}
\begin{proof}
For (i), if $\Gamma=\{1\}$, clearing the denominator of $c$ shows that $c\equiv1\pmod r$ for all but finitely many primes forces $c=1$. Otherwise, choose a $\Z$-basis of the free abelian group $\Gamma$. Since $\operatorname{End}_{\Q}(\mathbb G_m)=\Z$, this basis satisfies the hypothesis of~\cite[Theorem~8]{Perucca2010} for $G=\mathbb G_m/\Q$, yielding $c\in\Gamma$.

For (ii), suppose first that $\Gamma$ has positive rank, and enlarge $E$ to include all prime divisors of the numerators and denominators of $b$ and of a $\Z$-basis of $\Gamma$. If the conclusion failed, then $d_r$ would be coprime to $\ell$ for every $r\notin E$. By~\cite[Lemma~10]{Perucca2010} applied to $\mathbb G_m/\Q$, there would exist a nonzero integer $c$ (coprime to $\ell$) such that $b^c\in\Gamma$, contradicting the assumption that no positive power of $b$ lies in $\Gamma$ (taking $|c|$ if $c<0$).

If $\Gamma=\{1\}$, choose a prime $a$ dividing neither the numerator nor the denominator of $b$. No positive power of $b$ belongs to $\langle a\rangle$: if $b^m=a^j$ with $m>0$, unique factorization forces $j=0$ and hence $b^m=1$, contrary to hypothesis. Applying the positive-rank case to $\langle a\rangle$ gives a prime $r\notin E$ where $\overline b_r\langle\overline a_r\rangle$ has order divisible by $\ell$; since this order divides the order of $\overline b_r$ in $\mathbb F_r^\times$, the latter is also divisible by $\ell$.
\end{proof}

\medskip
\noindent\textbf{From this point onward, we assume that $T=\PP\setminus S$ is finite.} Consequently, $\Gamma_T$ is finitely generated, and every progression $a+d\Nzero$ with $\gcd(a,d)=1$ contains infinitely many primes in $S$.

The proof of cofinite rigidity has three main stages. We first obtain a multiplicative monoid automorphism $\mu$ by conjugating the power maps (\cref{thm:cofinite-power-conjugacy}); we then use $p$-adic closures and multiplicative orders in unit groups, whose cardinalities are given by Euler's totient function, to fix all infinite-height primes and prove that $\mu=\operatorname{id}$ (\cref{prop:cof-infinite-coordinates-fixed,prop:cof-mu-identity}). Finally, we rule out exchanges between finite-height primes with the same totient value (\cref{prop:cof-finite-coordinates-fixed}).

Our next goal is to show that $h$ conjugates the power maps $P_n(x)=x^n$ to power maps, extending~\cite[Theorem~5.1(4) and Lemma~5.15]{Banakh2018}. Note that for each $n\ge1$, $P_n\colon X_e\to X_e$ is continuous: if $x^n\in a+d\Nzero$, then $\gcd(x,d)=1$, so $x+d\Nzero$ is an admissible neighborhood of $x$ mapped into $x^n+d\Nzero\subseteq a+d\Nzero$ (for the Golomb, Kirch, and constant finite-height cases, see~\cite{Szczuka2013Darboux,Banakh2018,Szyszkowska2022b}).

\begin{theorem}
\label{thm:cofinite-power-conjugacy}
There exists a monoid automorphism $\mu\colon(\N,\cdot)\to(\N,\cdot)$ such that
\begin{equation}\label{eq:cofinite-power-conjugacy}
 h(x^n)=h(x)^{\mu(n)}\qquad(x,n\in\N).
\end{equation}
In particular, $h(1)=1$.
\end{theorem}
\begin{proof}
\noindent\textbf{Step 1.} We first express $h(a^n)$ as a power of $h(a)$ up to a factor in $\Gamma_T$. Fix $n\ge1$, put $\Phi=h\circ P_n\circ h^{-1}$, and take $a\in \mathcal Y_T$, $b=h(a)$, $c=h(a^n)$. Since $\Pi_S(a^n)=\Pi_S(a)$,~\eqref{eq:active-support-transport} gives $\Pi_S(b)=\Pi_S(c)=\sigma[\Pi_S(a)]\ne\varnothing$. We claim that
\begin{equation}\label{eq:localized-T-power-ray}
 c=u_a(n)b^{m_a(n)},\qquad u_a(n)\in\Gamma_T,\quad m_a(n)\in\Z.
\end{equation}
Suppose not. By \cref{lem:external-finite-rank}(i) applied to $\langle\Gamma_T,b\rangle$, there exists $r\in S$ with $r\nmid bc$ such that $\overline c_r\notin\langle\overline{\Gamma_T}_r,\overline b_r\rangle$. Since $\Phi(b)=c$, continuity yields an admissible progression $b+d\Nzero\subseteq C_e(b,r,1)\cap\Phi^{-1}[C_e(c,r,1)]$. Because $T$ is finite, Dirichlet's theorem gives a prime $s\in S\cap(b+d\Nzero)$. By~\eqref{eq:active-support-transport}, $\Pi_S(h^{-1}(s)^n)=\Pi_S(h^{-1}(s))=\{\sigma^{-1}(s)\}$, so $\Pi_S(\Phi(s))=\{s\}$, which means $\Phi(s)=vs^j$ for some $v\in\mathcal U_T$ and $j\ge1$. The congruences $s\equiv b\pmod r$ and $vs^j=\Phi(s)\equiv c\pmod r$ then imply $\overline c_r\in\langle\overline{\Gamma_T}_r,\overline b_r\rangle$, a contradiction. Comparing $q$-adic valuations in~\eqref{eq:localized-T-power-ray} at any $q\in\Pi_S(b)$ shows that $m_a(n)$ is unique and positive.

\smallskip\noindent\textbf{Step 2.} We show that $m_a(n)$ is independent of $a\in \mathcal Y_T$. Fix a prime $\ell$. Since $b=h(a)\in \mathcal Y_T$, no positive power of $b$ lies in $\Gamma_T$. By \cref{lem:external-finite-rank}(ii), there is a prime $r\in S$ (with $r\nmid bc$ and $r\notin T$) such that $\overline b_r\,\overline{\Gamma_T}_r$ has order divisible by $\ell$ in $\mathbb F_r^\times/\overline{\Gamma_T}_r$. By continuity, on a sufficiently small relative neighborhood of $a$ in $\mathcal Y_T$, the reductions of $h(x)$ and $h(x^n)$ modulo $r$ agree with those of $b$ and $c$, respectively. Comparing~\eqref{eq:localized-T-power-ray} at $x$ and at $a$ modulo $r$ yields
\[
 \overline b_r^{\,m_x(n)-m_a(n)}\in\overline{\Gamma_T}_r,
 \qquad\text{so}\qquad m_x(n)\equiv m_a(n)\pmod\ell.
\]
Thus $x\mapsto[m_x(n)]_\ell$ is locally constant on the connected space $\mathcal Y_T$ (\cref{lem:connected-active-layer}), hence constant. As $\ell$ was arbitrary, $m_a(n)$ does not depend on $a\in \mathcal Y_T$; denote its common value by $\mu(n)$.

\smallskip\noindent\textbf{Step 3.} We next show that $u_a(n)$ is independent of $a\in \mathcal Y_T$. Write $u_a(n)=\prod_{t\in T}t^{z_t(a,n)}$. For each $t\in T$ and any prime $\ell$, apply \cref{lem:external-finite-rank}(ii) to $t$ modulo $\Gamma_{T\setminus\{t\}}$ to obtain $r\in S$ (with $r\nmid bc$) at which the coset of $\overline t_r$ modulo $\overline{\Gamma_{T\setminus\{t\}}}_r$ has order divisible by $\ell$. On a relative neighborhood of $a$ in $\mathcal Y_T$ where the reductions of $h(x)$ and $h(x^n)$ modulo $r$ are constant, the reduction of $u_x(n)=h(x^n)h(x)^{-\mu(n)}$ modulo $r$ is also constant. Reducing $u_x(n)/u_a(n)$ modulo $r$ and quotienting by $\overline{\Gamma_{T\setminus\{t\}}}_r$ gives $z_t(x,n)\equiv z_t(a,n)\pmod\ell$. By connectedness of $\mathcal Y_T$ and variation of $\ell$, each $z_t(a,n)$ is independent of $a\in \mathcal Y_T$. Thus there exists $c_n\in\Gamma_T$ such that
\begin{equation}\label{eq:constant-T-factor-global}
 h(a^n)=c_n h(a)^{\mu(n)}\qquad(a\in \mathcal Y_T).
\end{equation}
Choosing $a\in \mathcal Y_T$ with $h(a)\in S$ via~\eqref{eq:T-unit-invariance}, any prime divisor of the denominator of $c_n$ lies in $T$ and cannot cancel with $h(a)^{\mu(n)}$; since $h(a^n)\in\N$, we must have $c_n\in\mathcal U_T$.

\smallskip\noindent\textbf{Step 4.} We extend~\eqref{eq:constant-T-factor-global} to all of $\N$ and show that $c_n=1$. Multiplication by $c_n\in\mathcal U_T$ is continuous on $X_e$ because $c_n$ is coprime to every admissible common difference. Thus both sides of~\eqref{eq:constant-T-factor-global} define continuous maps $X_e\to X_e$ that agree on the dense subset $\mathcal Y_T$; since $X_e$ is Hausdorff, they agree on all of $\N$. Setting $a=1$ gives $c_n=h(1)^{1-\mu(n)}$. If $\mu(2)=1$, then $h(a^2)=h(a)$ for all $a\in\N$, contradicting injectivity of $h$. Hence $\mu(2)>1$. Since $c_2=h(1)^{1-\mu(2)}$ is an integer, we must have $h(1)=1$, and therefore $c_n=1$ for every $n\ge1$.

\smallskip\noindent\textbf{Step 5.} Evaluating $h(x^{mn})=h((x^m)^n)$ at any $x$ with $h(x)>1$ gives $\mu(mn)=\mu(m)\mu(n)$. Applying Steps~1--4 to $h^{-1}$ produces an analogous multiplicative map $\nu$. For any $x>1$, the identities $x^n=h^{-1}(h(x^n))=x^{\nu(\mu(n))}$ and $x^n=h(h^{-1}(x^n))=x^{\mu(\nu(n))}$ show that $\nu\circ\mu=\operatorname{id}_{\N}$ and $\mu\circ\nu=\operatorname{id}_{\N}$, so $\nu=\mu^{-1}$. Thus $\mu$ is a monoid automorphism of $(\N,\cdot)$.
\end{proof}

To prove that $\sigma$ fixes every infinite-height prime, we examine the $p$-adic closures of sets of positive powers~\cite[Lemmas~8, 10, 12, 14, 16, and~20]{Banakh2021a}. Equip $U_p=\N\setminus p\N$ with the $p$-adic topology, and for $A\subseteq U_p$ write $\overline A^{\,p}=\operatorname{cl}_{U_p}(A)$. Define
\[
 \mathcal X_p=\left\{\overline{x^{\N}}^{\,p}:x\in U_p\setminus\{1\}\right\},
\]
ordered by \emph{reverse inclusion}: $X\le_p Y\iff Y\subseteq X$. Under this order, $\mathord{\uparrow}_p X=\{Y\in\mathcal X_p:Y\subseteq X\}$. We call $X\in\mathcal X_p$ an \emph{upper-chain element} if $\mathord{\uparrow}_p X$ is linearly ordered, and an \emph{antichain} is a subset of pairwise incomparable elements.

\begin{lemma}\label{lem:external-monogenic-structure}
Let $p$ be a prime.
\begin{enumerate}
\item If $p$ is odd, the upper-chain elements of $\mathcal X_p$ are precisely the sets $1+p^k\Nzero$ ($k\in\N$), which are well-ordered by $\le_p$:
\[
 1+p^k\Nzero\le_p 1+p^\ell\Nzero\quad\Longleftrightarrow\quad k\le\ell.
\]
\item The poset $\mathcal X_2$ contains an infinite antichain, whereas $\mathcal X_p$ contains no infinite antichain for any odd prime $p$.
\end{enumerate}
\end{lemma}
\begin{proof}
Part~(i) is~\cite[Lemma~14]{Banakh2021a}.

For~(ii), if $p$ is odd, \cite[Lemmas~10 and~12]{Banakh2021a} identify $\mathcal X_p$ with the divisibility poset $\mathcal D_p=\{d\in\N:d\mid p^n(p-1)\text{ for some }n\in\Nzero\}$. Every element of $\mathcal D_p$ has the form $p^j c$ with $c\mid p-1$. Since $p-1$ has only finitely many divisors, any infinite subset of $\mathcal D_p$ contains two elements $p^j c$ and $p^k c$ with $j<k$, which are comparable; hence $\mathcal X_p$ has no infinite antichain.

For $p=2$, \cite[Lemma~16(2)]{Banakh2021a} shows that
\[
 Y_n=\overline{(2^n-1)^{\N}}^{\,2}=(1+2^{n+1}\Nzero)\cup(2^n-1+2^{n+1}\Nzero)
 \qquad(n\ge2)
\]
belongs to $\mathcal X_2$. Every $y\in Y_n$ with $y\equiv3\pmod4$ satisfies $v_2(y+1)=n$. In particular, $2^n-1\in Y_n\setminus Y_m$ whenever $m\ne n$, so $\{Y_n:n\ge2\}$ is an infinite antichain in $\mathcal X_2$.
\end{proof}

Fix the monoid automorphism $\mu$ from \cref{thm:cofinite-power-conjugacy}. Since $\mu$ is bijective, $h[x^\N]=h(x)^\N$ for every $x\in\N$. For $p\in I$ and $q=\sigma(p)$, \cref{cor:cof-padic-reflection} therefore induces an order isomorphism
\begin{equation}\label{eq:cof-monogenic-transport}
 \mathcal X_p\longrightarrow\mathcal X_q,\qquad
 \overline{x^{\N}}^{\,p}\longmapsto\overline{h(x)^{\N}}^{\,q}.
\end{equation}
By \cref{lem:external-monogenic-structure}(ii), $\mathcal X_2$ cannot be order-isomorphic to $\mathcal X_q$ for any odd $q$, so $\sigma(2)=2$ whenever $2\in I$. For odd $p\in I$ and $q=\sigma(p)$, the upper-chain elements form an $\omega$-chain under $\le_p$, so any order isomorphism preserves each level:
\begin{equation}\label{eq:cof-level-transport}
 h[1+p^k\Nzero]=1+q^k\Nzero\qquad(k\ge1).
\end{equation}

\begin{lemma}
\label{lem:cof-infinite-order-transport}
For odd $p\in I$, $q=\sigma(p)$, and $p\nmid x$,
\begin{equation}\label{eq:cof-infinite-order}
 \ord_{q^k}(h(x))=\mu\bigl(\ord_{p^k}(x)\bigr)\qquad(k\ge1).
\end{equation}
Consequently,
\begin{equation}\label{eq:cof-p-minus-one-transport}
 \mu(p)=q,\qquad \mu(p-1)=q-1.
\end{equation}
\end{lemma}
\begin{proof}
Put $A=\ord_{p^k}(x)$ and $B=\ord_{q^k}(h(x))$. By~\eqref{eq:cofinite-power-conjugacy} and~\eqref{eq:cof-level-transport},
\[
 A\mid n\quad\Longleftrightarrow\quad x^n\in1+p^k\Nzero
 \quad\Longleftrightarrow\quad h(x)^{\mu(n)}\in1+q^k\Nzero
 \quad\Longleftrightarrow\quad B\mid\mu(n).
\]
Since $\mu$ is a multiplicative automorphism, it preserves divisibility, forcing $B=\mu(A)$.

Because $(\Z/p^k\Z)^\times$ and $(\Z/q^k\Z)^\times$ are cyclic~\cite[Lemma~4(1)]{Banakh2021a}, applying~\eqref{eq:cof-infinite-order} to a generator of $(\Z/p^k\Z)^\times$ gives $\mu(\varphi(p^k))\mid\varphi(q^k)$, while applying the same argument to $h^{-1}$ gives the reverse divisibility. Thus $\mu(\varphi(p^k))=\varphi(q^k)$, and taking $k=1,2$ yields~\eqref{eq:cof-p-minus-one-transport}.
\end{proof}

\begin{lemma}
\label{lem:cof-prime-adic-completion}
For odd $p\in I$ and $q=\sigma(p)$, $\mu$ is a homeomorphism from $\N$ with the $p$-adic topology to $\N$ with the $q$-adic topology. Extending $\mu$ multiplicatively to $\Q_{>0}^\times$, its restriction to the subgroup of rationals of $p$-adic valuation zero extends to a topological group isomorphism
\begin{equation}\label{eq:cof-unit-completion}
 \widehat{\mu}_{p,q}\colon\Z_p^\times\longrightarrow\Z_q^\times.
\end{equation}
\end{lemma}
\begin{proof}
\noindent\textbf{Step 1.} Put $a=1+p$ and $b=h(a)>1$. For $t\in\Z_p$, $a^t=\exp(t\log a)$ defines a homeomorphism $\Z_p\to1+p\Z_p$~\cite[Theorems~8.7 and~8.13]{ConradPadicSeries}, since $v_p(\log a)=1$. By~\eqref{eq:cofinite-power-conjugacy} and~\eqref{eq:cof-level-transport}, $\overline{b^\N}^{\,q}=1+q\Nzero$; since $1+q\Nzero$ is dense in $1+q\Z_q$, it follows that the positive powers of $b$ are dense in $1+q\Z_q$, which forces $v_q(b-1)=1$ (as $v_q(b-1)\ge2$ would confine all powers to $1+q^2\Z_q$). By the lifting-the-exponent lemma, for all $m\ne n$,
\[
 v_p(a^m-a^n)=1+v_p(m-n),\qquad
 v_q(b^m-b^n)=1+v_q(m-n),
\]
so $n\mapsto a^n$ ($\N\to a^\N$) and $n\mapsto b^n$ ($\N\to b^\N$) are homeomorphisms for the $p$-adic and $q$-adic topologies, respectively. By \cref{cor:cof-padic-reflection} and~\eqref{eq:cofinite-power-conjugacy}, $h$ restricts to a homeomorphism $a^\N\to b^\N$. Composing these maps shows that $\mu\colon\N\to\N$ is a homeomorphism from the $p$-adic to the $q$-adic topology.

\smallskip\noindent\textbf{Step 2.} Extend $\mu$ to a group automorphism $\Q_{>0}^\times\to\Q_{>0}^\times$, and put $\Gamma_p=\{u/v:u,v\in\N,\,p\nmid uv\}$. By~\eqref{eq:cof-p-minus-one-transport}, $\mu[\Gamma_p]=\Gamma_q$. For any $L\ge1$, continuity of $\mu$ at $1$ on $\N$ gives $K\ge1$ such that $n\equiv1\pmod{p^K}$ implies $\mu(n)\equiv1\pmod{q^L}$. If $u/v\in\Gamma_p$ satisfies $u\equiv v\pmod{p^K}$, choosing $c\in\N$ with $vc\equiv1\pmod{p^K}$ and applying the implication to $uc$ and $vc$ gives $\mu(u)/\mu(v)\equiv1\pmod{q^L}$. Thus $\mu|_{\Gamma_p}$ (and similarly its inverse) is continuous at the identity, hence uniformly continuous, and therefore extends to a topological group isomorphism $\widehat{\mu}_{p,q}\colon\Z_p^\times\to\Z_q^\times$.
\end{proof}

\begin{proposition}
\label{prop:cof-infinite-coordinates-fixed}
For every $p\in I$, $\sigma(p)=p$. For odd $p\in I$,
\begin{equation}\label{eq:cof-infinite-invariants}
 \mu(p)=p,\qquad \mu(p-1)=p-1.
\end{equation}
\end{proposition}
\begin{proof}
We already noted $\sigma(2)=2$ if $2\in I$. For odd $p\in I$ and $q=\sigma(p)$, the image $K=\widehat{\mu}_{p,q}[1+p\Z_p]$ is a compact (hence closed) infinite pro-$p$ subgroup of $\Z_q^\times$. Therefore $K\cap(1+q\Z_q)$ is a closed subgroup of both $K$ and $1+q\Z_q$, making it simultaneously a pro-$p$ and a pro-$q$ group. If $p\ne q$, every finite continuous quotient of $K\cap(1+q\Z_q)$ has both $p$-power and $q$-power order and is thus trivial; since finite continuous quotients separate points in a profinite group, $K\cap(1+q\Z_q)=\{1\}$. Reduction modulo $q$ would then embed the infinite group $K$ into the finite group $(\Z/q\Z)^\times$, a contradiction. Thus $p=q$, and~\eqref{eq:cof-p-minus-one-transport} gives~\eqref{eq:cof-infinite-invariants}.
\end{proof}

\begin{lemma}
\label{lem:cof-coordinate-order-invariants}
For $S_0=\{p\in S:p\ne2,\ \sigma(p)\ne2\}$, define
\begin{equation}\label{eq:Ap-definition}
 A_p=\begin{cases}
 p-1,&p\in I,\\
 \varphi(p^{e(p)}),&p\in F.
 \end{cases}
\end{equation}
Then
\begin{equation}\label{eq:Ap-fixed}
 \mu(A_p)=A_p\qquad(p\in S_0).
\end{equation}
\end{lemma}
\begin{proof}
For $p\in I$, this is~\eqref{eq:cof-infinite-invariants}. For $p\in F\cap S_0$, put $q=\sigma(p)$. Since $h(1)=1$, we have $\theta_p(1)=1$. For any $u\in(\Z/p^{e(p)}\Z)^\times$ represented by $x\in\N$, put $A=\ord_{p^{e(p)}}(u)$ and $B=\ord_{q^{e(q)}}(\theta_p(u))$. By~\eqref{eq:cof-terminal-transport} and~\eqref{eq:cofinite-power-conjugacy},
\[
 A\mid n\quad\Longleftrightarrow\quad x^n\equiv1\pmod{p^{e(p)}}
 \quad\Longleftrightarrow\quad h(x)^{\mu(n)}\equiv1\pmod{q^{e(q)}}
 \quad\Longleftrightarrow\quad B\mid\mu(n),
\]
so $B=\mu(A)$. Applying this to a generator of the cyclic group $(\Z/p^{e(p)}\Z)^\times$ (and similarly for $h^{-1}$) gives $\mu(\varphi(p^{e(p)}))=\varphi(q^{e(q)})=\varphi(p^{e(p)})$ by~\eqref{eq:cof-totient-cardinality}.
\end{proof}

\begin{proposition}
\label{prop:cof-mu-identity}
The monoid automorphism $\mu$ is the identity on $\N$.
\end{proposition}
\begin{proof}
If $\mu(2)=r\ne2$, then $r$ is odd; since $S_0$ is cofinite, Dirichlet's theorem gives $p\in S_0$ with $p\equiv2\pmod r$. Then $2\mid A_p$, so $r\mid\mu(A_p)=A_p$ by~\eqref{eq:Ap-fixed}. Yet every prime divisor of $A_p$ divides $p(p-1)$, whereas $r\nmid p(p-1)$, a contradiction. Thus $\mu(2)=2$.

Similarly, if $\mu(\ell)=r\ne\ell$ for an odd prime $\ell$, then $r$ is odd (since $\mu(2)=2$ and $\mu$ is bijective). The Chinese remainder theorem and Dirichlet's theorem supply $p\in S_0$ with $p\equiv1\pmod\ell$ and $p\equiv2\pmod r$. Then $\ell\mid A_p$ implies $r\mid A_p$, again contradicting $r\nmid p(p-1)$. Hence $\mu$ fixes every prime.
\end{proof}

Consequently,
\begin{equation}\label{eq:exact-power-compatibility}
 h(x^n)=h(x)^n\qquad(x,n\in\N).
\end{equation}

\begin{proposition}
\label{prop:cof-multiplicativity}
Every $h\in\Homeo(X_e)$ is completely multiplicative:
\begin{equation}\label{eq:h-multiplicative}
 h(xy)=h(x)h(y)\qquad(x,y\in\N).
\end{equation}
\end{proposition}
\begin{proof}
Suppose first that $I$ is infinite. Fix $x,y\in\N$ and choose an odd $p\in I$ with $p\nmid xy$. By \cref{prop:cof-infinite-coordinates-fixed} and \cref{cor:cof-padic-reflection}, $h|_{U_p}$ is a $p$-adic self-homeomorphism. Choose a primitive root $g>1$ modulo $p^2$; its powers are dense in $\Z_p^\times$~\cite[Lemma~9]{Banakh2021a}. Taking sequences $g^{a_j}\to x$ and $g^{b_j}\to y$ in $\Z_p$,~\eqref{eq:exact-power-compatibility} gives
\[
 h(g^{a_j+b_j})=h(g)^{a_j+b_j}=h(g^{a_j})h(g^{b_j}).
\]
Passing to the $p$-adic limit yields $h(xy)=h(x)h(y)$.

Suppose now that $I$ is finite, so $F$ is cofinite. For every $p\in F\cap S_0$,~\eqref{eq:cof-terminal-transport} and~\eqref{eq:exact-power-compatibility} give $\theta_p(u^n)=\theta_p(u)^n$. Since $\theta_p\colon(\Z/p^{e(p)}\Z)^\times\to(\Z/\sigma(p)^{e(\sigma(p))}\Z)^\times$ is a bijection between cyclic groups of the same order that preserves powers, it sends a generator to a generator and is therefore a group isomorphism. For every $p\in F\cap S_0$ with $p\nmid xy$ and $q=\sigma(p)$, we thus have $h(xy)\equiv h(x)h(y)\pmod{q^{e(q)}}$. Since there are infinitely many such primes $q$, the integer $h(xy)-h(x)h(y)$ must be zero.
\end{proof}

\begin{corollary}
\label{cor:cof-exact-prime-images}
Every $h\in\Homeo(X_e)$ is a monoid automorphism of $(\N,\cdot)$, and
\begin{equation}\label{eq:cof-exact-prime-images}
 h(p)=\sigma(p)\qquad(p\in S).
\end{equation}
\end{corollary}
\begin{proof}
Since $h$ is a multiplicative bijection fixing $1$, it permutes the primes. Equation~\eqref{eq:active-support-transport} then forces $h(p)=\sigma(p)$ for all $p\in S$.
\end{proof}

It remains to show that $\sigma$ fixes every finite-height prime $p\in F$. We already know that $\sigma$ preserves $\varphi(p^{e(p)})$. The next two elementary lemmas describe when two distinct prime powers can have the same totient value.

\begin{lemma}
\label{lem:cof-euler-collision}
If distinct primes $p<q$ and integers $a,b\ge1$ satisfy $\varphi(p^a)=\varphi(q^b)$, then
\begin{equation}\label{eq:euler-collision-shape}
 a\ge2,\qquad b=1,\qquad q-1=(p-1)p^{a-1}.
\end{equation}
In particular, any given totient value occurs at prime powers of at most two distinct prime bases.
\end{lemma}
\begin{proof}
In $p^{a-1}(p-1)=q^{b-1}(q-1)$, the prime $q>p$ cannot divide the left-hand side, so $b=1$. If $a=1$ as well, then $p-1=q-1$, contradicting $p<q$. This proves~\eqref{eq:euler-collision-shape}. If the same totient value occurred at three distinct primes, comparing the two larger primes with the smallest would force exponent $1$ at both larger primes, which is impossible since $\varphi(q_1)\ne\varphi(q_2)$ for distinct primes $q_1,q_2$.
\end{proof}

\begin{lemma}
\label{lem:cof-collision-stability}
If $q-1=(p-1)p^{a-1}$ for distinct primes $p<q$ and $a\ge2$, then for every prime $\ell\notin\{p,q\}$,
\[
 p\equiv1\pmod\ell\quad\Longleftrightarrow\quad q\equiv1\pmod\ell.
\]
\end{lemma}
\begin{proof}
Since $\ell\ne p$, $\ell\mid q-1$ if and only if $\ell\mid p-1$.
\end{proof}

\begin{proposition}
\label{prop:cof-finite-coordinates-fixed}
For every $p\in F$, $\sigma(p)=p$.
\end{proposition}
\begin{proof}
\noindent\textbf{Step 1.} Suppose for contradiction that $\sigma$ moves some prime in $F$. If $e(2)=1$, then $\sigma(2)=2$ by \cref{lem:cof-parity-repair}. On $F\cap S^\sharp$, $\sigma$ preserves $\varphi(p^{e(p)})$, and by \cref{lem:cof-euler-collision} each totient value occurs at no more than two prime bases. Hence any moved prime in $F$ must belong to a $2$-cycle of $\sigma$. Labeling this pair so that $p<q$, we have
\begin{equation}\label{eq:moved-collision-pair}
 \sigma(p)=q,\quad\sigma(q)=p,\quad
 e(p)=a\ge2,\quad e(q)=1,\quad q-1=(p-1)p^{a-1}.
\end{equation}

\smallskip\noindent\textbf{Step 2.} Since $q>p\ge2$, $q$ is odd, so $\gcd(2,q)=1$ and $\gcd(p^a,q)=1$. Because $T$ is finite, the Chinese remainder theorem and Dirichlet's theorem yield a prime $r\in S\setminus\{2,p,q\}$ satisfying
\begin{equation}\label{eq:third-prime-choice}
 r\equiv1\pmod{p^a},\qquad r\equiv2\pmod q.
\end{equation}
Since $h(1)=1$, the bijection $\theta_p$ in~\eqref{eq:cof-terminal-transport} sends the identity class modulo $p^a$ to the identity class modulo $q$. Together with \cref{cor:cof-exact-prime-images}, this gives
\begin{equation}\label{eq:sigma-r-one-mod-q}
 \sigma(r)=h(r)\equiv1\pmod q.
\end{equation}
Comparing~\eqref{eq:third-prime-choice} with~\eqref{eq:sigma-r-one-mod-q} shows that $\sigma(r)\ne r$.

\smallskip\noindent\textbf{Step 3.} By \cref{prop:cof-infinite-coordinates-fixed}, all infinite-height primes are fixed by $\sigma$, so $r$ and $\sigma(r)$ are distinct primes in $F$ with equal totient values. Neither equals $p$ or $q$, since $\{p,q\}$ already forms a $2$-cycle of $\sigma$. Applying \cref{lem:cof-collision-stability} to the pair $\{r,\sigma(r)\}$ with $\ell=q$ yields $r\equiv1\pmod q$ from~\eqref{eq:sigma-r-one-mod-q}, contradicting~\eqref{eq:third-prime-choice}. Thus $\sigma(p)=p$ for all $p\in F$.
\end{proof}

\begin{lemma}
\label{lem:active-primes-dense}
The set $S$ of active primes is dense in $X_e$.
\end{lemma}
\begin{proof}
Every basic progression $a+d\Nzero$ with $\gcd(a,d)=1$ contains infinitely many primes by \cref{lem:dirichlet}, and only finitely many primes lie outside $S$. (For the Golomb case, see~\cite[Theorem~6]{Golomb1959} and~\cite[Theorem~3.1]{Banakh2018}.)
\end{proof}

\begin{theorem}
\label{thm:main-cofinite-arbitrary}
If $\PP\setminus S_e$ is finite, then $X_e$ is rigid.
\end{theorem}
\begin{proof}
Let $h\in\Homeo(X_e)$. By \cref{prop:cof-infinite-coordinates-fixed,prop:cof-finite-coordinates-fixed} and \cref{cor:cof-exact-prime-images}, $h(p)=\sigma(p)=p$ for every $p\in S_e$. Since $S_e$ is dense in $X_e$ (\cref{lem:active-primes-dense}) and $X_e$ is Hausdorff (\cref{prop:general-specialization}), $h$ is the identity on $\N$.
\end{proof}

\begin{corollary}
\label{cor:full-support-interpolation}
If $e(p)\ge1$ for every prime $p$, then $X_e$ is rigid.
\end{corollary}
\begin{proof}
Apply \cref{thm:main-cofinite-arbitrary} with $T=\varnothing$.
\end{proof}

\section{Nonrigidity on a support of density one}
\label{sec:general-infinite-supports}

By \cref{thm:main-cofinite-arbitrary}, $X_e$ is rigid whenever all but finitely many primes are active. We now show that deleting a set of relative density zero can give a space with a fixed-point-free involution. Throughout this section, $q_j$ denotes the $j$-th prime and $\pi(x)$ the number of primes at most $x$. All densities of sets of primes are relative to $\pi(x)$.

\begin{theorem}\label{thm:sparse-involution}
There exist a coinfinite set $G\subseteq\PP$ of natural density one among the primes and a fixed-point-free involution $h$ of $\N$ with $h(1)=2$ such that the following holds. For every infinite set $S\subseteq G$, the map $h$ is a self-homeomorphism of $X_e$ for the height function equal to $1$ on $S$ and $0$ elsewhere.
\end{theorem}
\begin{proof}
We construct $h$ by pairing positive integers. We include $0$ as an auxiliary fixed point so that the induced permutations of residue classes fix the zero class.

\smallskip\noindent\textbf{Step 1. Extending a finite pairing.}
Start with $A_0=\{0,1,2\}$, let $f_0$ fix $0$ and exchange $1$ and $2$, and set $E_0=\{2\}$. For every $p\notin E_0$, the elements of $A_0$ have distinct residues modulo $p$. At stage $s\ge1$, let $A=A_{s-1}$ be the current finite set and $f=f_{s-1}$ its involution. Our induction hypothesis is that for every prime $p$ outside $E_0\cup\cdots\cup E_{s-1}$,
\begin{equation}\label{eq:sparse-involution-invariant}
 a\equiv b\pmod p\ \Longrightarrow\ f(a)\equiv f(b)\pmod p
 \qquad(a,b\in A).
\end{equation}
Thus $f$ induces an involution on the residue classes represented by $A$, fixing the class of $0$.

Let $n=\min(\N\setminus A)$, $U=A\cup\{n\}$ and $k=|U|$. Choose a prime $B>\max U$ so large that, with $N=\pi(B)$ and $M=\lfloor B/(2k)\rfloor$,
\[
 M>N,\qquad \frac{k(N+2)}{M}\le2^{-s}.
\]
This is possible since $\pi(B)/B\to0$ and $k$ is fixed at this stage. Put $Q=\prod_{p\le B}p$. For each prime $p\le B$ outside the previous exceptional sets, extend the induced involution to $\theta_p$ on $\mathbb F_p$ by fixing the unrepresented classes. By the Chinese remainder theorem we may choose $1\le r\le Q$ satisfying
\begin{equation}\label{eq:sparse-new-partner}
 [r]_p=\theta_p([n]_p)
\end{equation}
at all these primes. Every candidate $m=r+tQ$, $1\le t\le B$, exceeds $B$ and satisfies the same congruences.

Choose a candidate for which $m\not\equiv a\pmod p$ whenever $a\in U$ and $B<p\le q_M$. For a given pair $(a,p)$, at most one value of $t$ is forbidden, since $p\nmid Q$ and $p>B$. There are at most $kM\le B/2$ such pairs, so at least one candidate remains. Set $A_s=U\cup\{m\}$ and extend $f$ to $f_s$ by exchanging $n$ and $m$.

\smallskip\noindent\textbf{Step 2. Bounding the exceptional primes.}
Let
\[
 E_s=\left\{p\in\PP:p>q_M,\quad p\mid\prod_{a\in U}(m-a)\right\}.
\]
At each previously retained prime $p\le B$, the new pair respects $\theta_p$ in both directions by~\eqref{eq:sparse-new-partner} and $\theta_p^2=\operatorname{id}$. For $p>B$, the elements of $U$ have distinct residues modulo $p$ (since $B>\max U$). By the choice of $m$ for $p\le q_M$ and the definition of $E_s$ for $p>q_M$, these residues are distinct from $[m]_p$ whenever $p\notin E_s$. Thus~\eqref{eq:sparse-involution-invariant} holds on $A_s$ outside $E_0\cup\cdots\cup E_s$.

Since $B\ge3$ and $Q\le B^N$, we have $m\le(B+1)Q\le B^{N+2}$, so the product defining $E_s$ is at most $B^{k(N+2)}$. Every prime in $E_s$ exceeds $B$, so $|E_s|\le k(N+2)$. For $x<q_M$ we have $E_s\cap[1,x]=\varnothing$, while for $x\ge q_M$ we have $\pi(x)\ge M$. Consequently
\[
 \frac{|E_s\cap[1,x]|}{\pi(x)}\le\frac{k(N+2)}{M}\le 2^{-s}
 \qquad(x\ge2).
\]

\smallskip\noindent\textbf{Step 3. Passing to the full involution.}
Put $E=\bigcup_{s\ge0}E_s$ and $G=\PP\setminus E$. For each $L\ge1$, the union of the first $L+1$ exceptional sets is finite. The uniform estimate of Step~2 therefore gives
\[
 \limsup_{x\to\infty}\frac{|E\cap[1,x]|}{\pi(x)}
 \le \sum_{s>L}2^{-s}=2^{-L}.
\]
Letting $L\to\infty$ shows that $E$ has density zero, so $G$ has density one among the primes. Each stage pairs the least unused positive integer, so the $f_s$ exhaust $\Nzero$. Their union $H$ fixes $0$, and its restriction $h$ to $\N$ is a fixed-point-free involution with $h(1)=2$. For $p\in G$, any two integers in $\Nzero$ occur together in some $A_s$. By~\eqref{eq:sparse-involution-invariant}, $H$ therefore induces an involution $\theta_p$ of $\mathbb F_p$ fixing $0$, with
\begin{equation}\label{eq:sparse-global-congruence}
 [h(n)]_p=\theta_p([n]_p)\qquad(n\in\N).
\end{equation}

Let $S\subseteq G$ be infinite and let $e$ be as in the statement. By \cref{prop:inf-terminal-basis}, the cylinders $C_e(u,p,1)$ with $p\in S$ and $u\in\mathbb F_p^\times$ form a subbasis for $\taue$. By~\eqref{eq:sparse-global-congruence}, the preimage of $C_e(u,p,1)$ under $h$ is $C_e(\theta_p(u),p,1)$, which is open since $\theta_p(u)\ne0$. Thus $h$ is continuous, and $h^{-1}=h$ gives continuity of its inverse. Taking $S=G$, we see from \cref{thm:main-cofinite-arbitrary} that $\PP\setminus G$ must be infinite.
\end{proof}

\begin{remark}\label{ex:sparse-nonrigid}
For instance, $S=G\cap\{q_{2j}:j\ge1\}$ has natural density $1/2$ among the primes. Indeed, the even-indexed primes contribute $\lfloor\pi(x)/2\rfloor$ terms up to $x$, and intersecting with $G$ removes at most $o(\pi(x))$ of them. The same involution is a self-homeomorphism for this support.
\end{remark}

It remains open to characterize the infinite, coinfinite sets $S$ for which $X_e$ is rigid when $e=1$ on $S$ and $e=0$ elsewhere.


\section*{Declaration of generative AI and AI-assisted technologies in the manuscript preparation process}

During the preparation of this work, the author used ChatGPT and Gemini to help review mathematical arguments, search for counterexamples, verify citations, and edit the manuscript. After using these tools, the author reviewed and edited the content as needed and takes full responsibility for the content of the article.

\bibliographystyle{abbrv}
\bibliography{lili}
\end{document}